\documentclass[11pt, leqno]{article}

\usepackage{amsthm}
\usepackage{amsmath}
\usepackage{amsfonts}
\usepackage{amssymb}

\usepackage{amsmath,amsthm,amscd,amssymb}
\usepackage{latexsym}
\usepackage{graphicx}
\usepackage[all]{xy}

\newtheorem{thm}{\textbf{Theorem}}

\newtheorem{lem}{\textbf{Lemma}}
\newtheorem{df}{\textbf{Definition}}
\newtheorem{cor}{\textbf{Corollary}}
\newtheorem{prop}{\textbf{Proposition}}

\newcommand{\A}{\Bbb{A}}
\newcommand{\R}{\Bbb{R}}

\newcommand{\C}{\Bbb{C}}

\newcommand{\Q}{\Bbb{Q}}

\newcommand{\OO}{\mathcal{O}}
\newcommand{\F}{\Bbb{F}}

\newcommand{\Z}{\Bbb{Z}}

\newcommand{\T}{\Bbb{T}}
\newcommand{\p}{\frak{p}}
\newcommand{\q}{\frak{q}}

\newcommand{\h}{\frak{t}}
\newcommand{\m}{\frak{m}}

\newcommand{\diag}{\text{diag}}

\newcommand{\Ind}{\text{Ind}}

\newcommand{\Ann}{\text{Ann}}

\newcommand{\rk}{\text{rk}}

\newcommand{\End}{\text{End}}

\newcommand{\Hom}{\text{Hom}}

\newcommand{\la}{\langle}
\newcommand{\ra}{\rangle}

\newcommand{\im}{\text{im}}

\newcommand{\Fr}{\text{Fr}}
\newcommand{\GL}{\text{GL}}
\newcommand{\GSp}{\text{GSp}}

\newcommand{\GSpin}{\text{GSpin}}
\newcommand{\Sp}{\text{Sp}}

\newcommand{\SO}{\text{SO}}

\newcommand{\PGSp}{\text{PGSp}}

\newcommand{\SL}{\text{SL}}

\newcommand{\U}{\text{U}}

\newcommand{\e}{\vspace{1pt}}
\newcommand{\y}{\hspace{6pt}}
\newcommand{\et}{\vspace{6pt}}

\begin{document}

\title{\bf{A Generalization of Level-Raising \\Congruences for Algebraic Modular Forms
}}

\author{\sc{Claus Mazanti Sorensen} \\ Department of Mathematics
\\California Institute of Technology \\
{\texttt{claus@caltech.edu}}}

\date{}

\maketitle

\begin{abstract}
In this paper we start by extending the results of K. Ribet and R.
Taylor on level-raising for algebraic modular forms on
$D^{\times}$, where $D$ is a definite quaternion algebra over a
totally real field $F$. We do this for automorphic representations
$\pi$ of an arbitrary reductive group $G$ over $F$ which is
compact at infinity. If $\lambda$ is a finite place of $\bar{\Q}$,
and $w$ is a place where $\pi_w$ is unramified and $\pi_w \equiv
{\bf{1}}$ (mod $\lambda$), then under some mild additional
assumptions we prove the existence of a $\tilde{\pi} \equiv \pi$
(mod $\lambda$) such that $\tilde{\pi}_w$ has more parahoric fixed
vectors than $\pi_w$. In the case where $G_w$ has semisimple rank
one, we recover results due to L. Clozel and J. Bellaiche
according to which $\tilde{\pi}_w$ is ramified. To provide
applications of our main theorem we consider two examples over
$\Q$ of rank greater than one. In the first example we take $G$ to
be a unitary group in three variables, and in the second we take
$G$ to be an inner form of $\GSp(4)$. In both cases, we obtain
precise satisfiable conditions on a split prime $w$ guaranteeing
the existence of a $\tilde{\pi} \equiv \pi$ (mod $\lambda$) such
that the component $\tilde{\pi}_w$ is generic and Iwahori
spherical. For symplectic $G$, to conclude that $\tilde{\pi}_w$ is
generic, we use computations of R. Schmidt. In particular, if
$\pi$ is of Saito-Kurokawa type, it is congruent to a
$\tilde{\pi}$ which is not of Saito-Kurokawa type.

\end{abstract}

\raggedright

\section*{Introduction}

In this paper, we will prove a generalization of the following
theorem of Ribet [Rib]:

\begin{thm}
Let $f \in S_2(\Gamma_0(N))$ be an eigenform, and let
$\lambda|\ell$ be a finite place of $\bar{\Q}$ such that $\ell
\geq 5$ and $f$ is not congruent to an Eisenstein series modulo
$\lambda$. If $q \nmid N\ell$ is a prime number such that $\ell
\nmid 1+q$ and the following condition is satisfied,
$$
\text{$a_q(f)^2 \equiv (1+q)^2$ (mod $\lambda$),}
$$
then there exists a $q$-new eigenform $\tilde{f} \in
S_2(\Gamma_0(Nq))$ congruent to $f$ modulo $\lambda$.
\end{thm}

Two eigenforms $f$ and $\tilde{f}$ are said to be congruent modulo
$\lambda$ if their Hecke eigenvalues are congruent for almost all
primes, that is, if $a_p(f)\equiv a_p(\tilde{f})$ (mod $\lambda$)
for almost all $p$. The proof of this theorem can be reduced, via
the Jacquet-Langlands correspondence, to the corresponding
statement for $D^{\times}$ where $D$ is a definite quaternion
algebra over $\Q$.

$\e$

Our goal in this paper is to prove that an automorphic form of
Saito-Kurokawa type is congruent to an automorphic form which is
$not$ of Saito-Kurokawa type. Since functorialty is not yet
available, we are considering an inner form $G$ of $\PGSp(4)/\Q$
such that $G(\R)$ is compact. By a form on $G\simeq \SO(5)$ of
Saito-Kurokawa type we mean a theta lift from
$\widetilde{\SL}(2)$. We achieve this goal as a result of Theorem
7 in section 8.3 below.

$\e$

We apply some of the ideas and methods of [Ta1] and [Ta2]. The
level-raising part of Taylor's proof carries over to a much more
general setup. Namely, the following: We let $F$ denote a totally
real number field with adeles $\A=F_{\infty}\times \A^{\infty}$,
and let $G$ be a connected reductive $F$-group such that
$G_{\infty}^1=G_{\infty}\cap G(\A)^1$ is compact and
$G^{\text{der}}$ is simple and simply connected. When $F \neq \Q$,
this just means that $G_{\infty}$ is compact. However, when $F=\Q$
and $Z_G$ is split, it suffices that $G_{\infty}^{\text{der}}$ is
compact. There are plenty of such groups. In fact, any split
simple $F$-group not of type $A_n$ ($n \geq 2$), $D_{2n+1}$ or
$E_6$ has infinitely many inner forms which are compact at
infinity (and quasi-split at all but at most one finite place).
Throughout, we fix a Haar measure $\mu=\otimes \mu_v$ on
$G(\A^{\infty})$. It is convenient to state our results using the
following notion of congruence. As $K$ varies over the compact
open subgroups of $G(\A^{\infty})$, the centers
$Z(\mathcal{H}_{K,\Z})$ of the Hecke algebras form an inverse
system. To an automorphic representation $\pi$ of $G(\A)$ we
associate the character
$$
\eta_{\pi}:\underleftarrow \lim Z(\mathcal{H}_{K,\Z}) \rightarrow
\C
$$
such that $\eta_{\pi}=\eta_{\pi^K} \circ \text{pr}_K$ for every
compact open subgroup $K$ such that $\pi^K \neq 0$. If $\lambda$
is a finite place of $\bar{\Q}$, we say that $\tilde{\pi}$ and
$\pi$ are congruent modulo $\lambda$ if their characters are. We
write $\tilde{\pi}\equiv \pi$ (mod $\lambda$). A similar notion
makes sense locally, and then $\tilde{\pi}\equiv \pi$ (mod
$\lambda$) if and only if $\tilde{\pi}_v\equiv \pi_v$ (mod
$\lambda$)  for all finite $v$. Moreover, when both
$\tilde{\pi}_v$ and $\pi_v$ are unramified, $\tilde{\pi}_v\equiv
\pi_v$ (mod $\lambda$) simply means the Satake parameters are
congruent. Before we can state the main theorem, we need the
following definition.

\begin{df}
Let $\pi$ be an automorphic representation of $G(\A)$ such that
$\pi_{\infty}={\bf{1}}$. We say that $\pi$ is abelian modulo
$\lambda$, a finite place of $\bar{\Q}$, if there exists an
automorphic character $\chi$ of $G(\A)$ with
$\chi_{\infty}={\bf{1}}$ such that $\pi \equiv \chi$ (mod
$\lambda$).
\end{df}

This is the analogue of the notion Eisenstein modulo $\lambda$ in
[Clo, p. 1269]. Since $G^{\text{der}}$ is anisotropic in our
setup, there are no cusps and we prefer the terminology abelian
modulo $\lambda$. The following theorem is in some sense the main
result of this paper.

\begin{thm}
Let $\pi=\otimes \pi_v$ be an automorphic representation of
$G(\A)$ such that $\pi_{\infty}={\bf{1}}$, and let $\lambda|\ell$
be a finite place of $\bar{\Q}$ such that $\pi$ is non-abelian
modulo $\lambda$. Suppose $w$ is a finite place of $F$ where
$\pi_w$ is unramified and
$$
\text{$\pi_w \equiv {\bf{1}}$ (mod $\lambda$).}
$$
Let $K_w \subset G_w$ be a hyperspecial subgroup and let
$J_w=K_w\cap K_w'$ be a parahoric subgroup, where $K_w'$ is
another maximal compact subgroup. Suppose $\ell \nmid [K_w':J_w]$.
Then there exists an automorphic representation
$\tilde{\pi}=\otimes \tilde{\pi}_v$ of $G(\A)$ with
$\tilde{\pi}_{\infty}={\bf{1}}$,

\begin{itemize}

\item $\tilde{\pi}_w^{J_w}\neq \tilde{\pi}_w^{K_w}+\tilde{\pi}_w^{K_w'}$, \\
\item $\tilde{\pi}\equiv \pi$ (mod $\lambda$).

\end{itemize}

\end{thm}

This theorem has no content unless
$\pi_w^{J_w}=\pi_w^{K_w}+\pi_w^{K_w'}$. There is a more precise
version later in this paper. If $G_w^{\text{der}}$ has rank one,
$J_w$ is an Iwahori subgroup and one can conclude that
$\tilde{\pi}_w^{K_w}=0$ but $\tilde{\pi}_w^{J_w}\neq 0$. This was
first proved by Bellaiche in his thesis [Bel], using the ideas of
Clozel [Clo]. By a theorem of Serre, [Ser], the eigensystem of a
modular form mod $\ell$ comes from an algebraic modular form mod
$\ell$ on $D^{\times}$, where $D/\Q$ now is the quaternion algebra
with ramification locus $\{\infty,\ell\}$. Combining this result
with the Jacquet-Langlands correspondence yields the result of
Ribet after stripping powers of $\ell$ from the level.

$\e$

There is another proof of Ribet's theorem relying on the so-called
Ihara lemma. It states that for $q \nmid N\ell$, the degeneracy
maps $X_0(Nq) \rightrightarrows X_0(N)$ induce an injection
$$
H^1(X_0(N),\Z_{\ell})^{\oplus 2} \rightarrow
H^1(X_0(Nq),\Z_{\ell})
$$
with torsion-free cokernel. The proof of this lemma reduces to the
congruence subgroup property of the group $\SL_2(\Z[1/q])$. In our
case we are looking at functions on a finite set, and the analogue
of the Ihara lemma can be proved by imitating the combinatorial
argument of Taylor [Ta1, p. 274] in the diagonal weight $2$ case.
See section 4.3 below.

$\e$

We mention a few applications of our main theorem. First, let
$E/\Q$ be an imaginary quadratic extension and let $G^*=\U(2,1)$
be the quasi-split unitary $\Q$-group in $3$ variables split over
$E$. Let $G=\U(3)$ be an inner form of $G^*$ such that
$G_{\infty}$ is compact. For primes $q$ inert in $E$, the
semisimple rank of $G(\Q_q)$ is one and we recover the result of
Clozel [Clo]. In the split case we obtain the following as a
corollary:

\begin{thm}
Let $\pi=\otimes \pi_p$ be an automorphic representation of
$G(\A)$ with $\pi_{\infty}={\bf{1}}$, and let $\lambda|\ell$ be a
finite place of $\bar{\Q}$ such that $\pi$ is non-abelian modulo
$\lambda$. Suppose $q \neq \ell$ is a prime, split in $E$, such
that $\pi_q$ is unramified and $\ell \nmid 1+q+q^2$. If moreover,
for $\q|q$,
$$
\text{$\h_{\pi,\q}\equiv
\begin{pmatrix}q & &
\\ & 1 & \\ & & q^{-1}
\end{pmatrix}$ (mod $\lambda$),}
$$
then there exists an automorphic representation
$\tilde{\pi}=\otimes \tilde{\pi}_p$ of $G(\A)$ with
$\tilde{\pi}_{\infty}={\bf{1}}$,

\begin{itemize}

\item $\tilde{\pi}_q$ is generic and $\tilde{\pi}_q^{J_q}\neq 0$,
where $J_q$ is
any maximal proper parahoric, \\
\item $\tilde{\pi}\equiv \pi$ (mod $\lambda$).

\end{itemize}

\end{thm}

We cannot prove by our methods that $\tilde{\pi}_q$ is ramified.
On the other hand, Bellaiche has a result in his thesis in the
split case, [Bel, p. 218], proving that $\tilde{\pi}_q$ is
ramified under the additional assumption that $\pi$ occurs with
multiplicity one (and discarding finitely many primes $\ell$). We
classify the Iwahori-spherical representations of $\GL(3)$ and
compute the dimensions of their parahoric fixed spaces. This
allows us to conclude that $\tilde{\pi}_q$ is either a full
unramified principal series or induced from a Steinberg
representation.

$\e$

It seems very likely that our method and corollary can be extended
to allow $\pi_{\infty}\neq {\bf{1}}$, but we have chosen not to do
it here for the sake of brevity. In that case it would follow that
if $\pi$ is endoscopic abelian (that is, nearly equivalent to a
weak transfer of a character of an endoscopic group), then it is
congruent to a $\tilde{\pi}$ which is not endoscopic abelian. This
is true even for $\U(n)$, for all $n \geq 2$. For $n=3$ this
phenomenon has been applied to the Bloch-Kato conjecture for
certain Hecke characters of $E$ by Bellaiche [Bel].

$\e$

In our second application, we let $G$ be an inner form of
$\GSp(4)$ such that $G^{\text{der}}(\R)$ is compact. Concretely,
$G=\GSpin(f)$ for some definite quadratic form $f$ in $5$
variables over $\Q$. In this situation, our main theorem yields
the following:

\begin{thm}
Let $\pi=\otimes \pi_p$ be an automorphic representation of
$G(\A)$ with $\pi_{\infty}={\bf{1}}$, and let $\lambda|\ell$ be a
finite place of $\bar{\Q}$ such that $\pi$ is non-abelian modulo
$\lambda$. Suppose $q \neq \ell$ is a prime such that $\pi_q$ is
unramified. If the Hecke matrix satisfies
$$
\text{$\h_{\pi_q\otimes |\nu|^{-3/2}}\equiv
\begin{pmatrix}1 & & &
\\ & q &  & \\ & & q^2 & \\ & & & q^3
\end{pmatrix}$ (mod $\lambda$),}
$$
then there exists an automorphic representation
$\tilde{\pi}=\otimes \tilde{\pi}_p$ of $G(\A)$ with
$\tilde{\pi}_{\infty}={\bf{1}}$,

\begin{itemize}

\item $\tilde{\pi}_q$ is generic and $\tilde{\pi}_q^{J_q}\neq 0$,
where $J_q$ is
the Klingen parahoric,\\
\item $\tilde{\pi}\equiv \pi$ (mod $\lambda$).

\end{itemize}

\end{thm}

By the Klingen parahoric, we mean the inverse image of the
standard Klingen parabolic in $\GSp(4,\F_q)$ under the reduction
map. Briefly, the proof relies on the computations of Ralf Schmidt
[Sch]. If $m(\pi)=1$, Bellaiche's methods seem to apply and one
can probably show that $\tilde{\pi}_q$ is induced from a twisted
Steinberg representation on the standard Klingen-Levi. It is known
that Saito-Kurokawa lifts (that is, theta lifts from
$\widetilde{\SL}(2)$) are locally non-generic. Therefore, if $\pi$
is of Saito-Kurokawa type, it is congruent to a $\tilde{\pi}$
which is not of Saito-Kurokawa type. Our interest in it stems from
our desire to apply it to the Bloch-Kato conjecture for the
motives attached to classical modular forms, and we plan to study
this in a sequel paper. In particular, we hope to establish a mod
$\ell$ analogue of a result of Skinner and Urban [SU], which is
valid for $all$ (not necessarily ordinary) modular forms of
classical weight at least $4$.

$\e$

This work forms part of my doctoral dissertation at the California
Institute of Technology. I would like to acknowledge the impact of
the ideas of Ribet, Taylor, Clozel and Bellaiche.

\section{The Abstract Setup and Taylor's Lemma}

\subsection{The Abstract Setup}

In this section, we fix a subring $\OO \subset \C$ and denote by
$L \subset \C$ its field of fractions. Let $H$ be a commutative
$\C$-algebra. We do not require $H$ to be of finite dimension.
However, we assume $H$ comes equipped with an involution $\phi
\mapsto \phi^{\vee}$. Here, by involution we mean a $\C$-linear
automorphism of order two. Moreover, we fix an $\OO$-order
$H_{\OO}\subset H$ preserved by $\vee$. Then we look at a triple
$(V, \la -,- \ra_V, V_{\OO})$ consisting of the following data:

\begin{itemize}
\item $V$ is a finite-dimensional $\C$-space with an action
$r_V:H \rightarrow \End_{\C}(V)$, \\
\item $\la -,- \ra_V$ is a non-degenerate, symmetric,
$\C$-bilinear form
$V \times V \rightarrow \C$, \\
\item $V_{\OO}\subset V$ is an $\OO$-lattice (that is, the
$\OO$-span of a $\C$-basis).
\end{itemize}

We impose the following compatibility conditions on these data:

\begin{itemize}
\item $r_V(\phi^{\vee})$ is the adjoint of $r_V(\phi)$ with respect to $\la -,- \ra_V$,  \\
\item $V_{\OO}\subset V$ is preserved by the order $H_{\OO}\subset
H$, \\
\item $V_{\OO}/V_{\OO}\cap V_{\OO}^{\vee}$ and
$V_{\OO}^{\vee}/V_{\OO}\cap V_{\OO}^{\vee}$ are torsion
$\OO$-modules.
\end{itemize}

Here $V_{\OO}^{\vee}=\{v \in V: \la v,V_{\OO}\ra_V \subset \OO \}$
is the dual lattice of $V_{\OO}$ in $V$. Choose nonzero
annihilators $A_V $ and $B_V$ in $\OO$ of the torsion modules
above, that is, such that
$$
\text{$A_V \la V_{\OO},V_{\OO}\ra_V \subset \OO$ $\y$ and $\y$
$\la v,V_{\OO}\ra_V \subset \OO \Rightarrow B_Vv \in V_{\OO}$.}
$$
Now let $(U, \la -,- \ra_U, U_{\OO})$ be another such triple and
choose annihilators $A_U$ and $B_U$ for it as above. Suppose we
are given a map $\delta: U \rightarrow V$, which is $H$-linear,
and in addition has the following properties:

\begin{itemize}
\item $U=\ker(\delta)\oplus\ker(\delta)^{\bot}$,\\
\item $V=\im(\delta)\oplus\im(\delta)^{\bot}$,\\
\item $\delta(U_{\OO}) \subset V_{\OO} \cap \delta(U)$, and the
quotient is killed by $C \in \OO-\{0\}$.
\end{itemize}

We consider its adjoint map $\delta^{\vee}:V \rightarrow U$
defined in the obvious way. Let $V^{\text{old}}=\im(\delta)$ and
$V^{\text{new}}=\im(\delta)^{\bot}$. These are $H$-stable
subspaces of $V$, and by assumption we have an orthogonal
decomposition $V=V^{\text{old}} \oplus V^{\text{new}}$.

\begin{df}
Let $V_{\mathcal{O}}^{\text{old}}=V_{\mathcal{O}}\cap
V^{\text{old}}$ and
$V_{\mathcal{O}}^{\text{new}}=V_{\mathcal{O}}\cap V^{\text{new}}$.
\end{df}

These $H_{\mathcal{O}}$-stable submodules of $V_{\mathcal{O}}$
span $V^{\text{old}}$ and $V^{\text{new}}$ respectively. They are
orthogonal, but their sum is not always all of $V_{\mathcal{O}}$.
Note that $\delta(U_{\mathcal{O}}) \subset
V_{\mathcal{O}}^{\text{old}}$ and $C V_{\mathcal{O}}^{\text{old}}
\subset \delta(U_{\mathcal{O}})$ by assumption. Now we look at the
quotients of $\T_{\mathcal{O}}$, the image of $H_{\mathcal{O}}$ in
$\End_{\mathcal{O}}(V_{\mathcal{O}})$, cut out by these
submodules:
$$
\text{$\T_{\mathcal{O}}^{\text{old}} \subset
\End_{\mathcal{O}}(V_{\mathcal{O}}^{\text{old}})$ $\y$ and $\y$
$\T_{\mathcal{O}}^{\text{new}} \subset
\End_{\mathcal{O}}(V_{\mathcal{O}}^{\text{new}})$}
$$
denote the images of $H_{\mathcal{O}}$. Clearly we have natural
surjective maps $\T_{\mathcal{O}} \twoheadrightarrow
\T_{\mathcal{O}}^{\text{old}}$ and $\T_{\mathcal{O}}
\twoheadrightarrow \T_{\mathcal{O}}^{\text{new}}$ given by
restriction, and $\T_{\mathcal{O}}$ acts faithfully on
$V_{\mathcal{O}}$.

\subsection{An Extension of Taylor's Lemma}

Note that $\T_{\mathcal{O}}$ acts naturally on
$U_{\mathcal{O}}'=U_{\mathcal{O}}\cap \ker(\delta)^{\bot}$.
Moreover, one can easily check that the action factors through
$\T_{\mathcal{O}}^{\text{old}}$. By a congruence module we mean a
$\T_{\mathcal{O}}$-module, such that the action factors through
both $\T_{\mathcal{O}}^{\text{old}}$ and
$\T_{\mathcal{O}}^{\text{new}}$. The following lemma was stated
for $\mathcal{O}=\Z$, trivial annihilators, and injective $\delta$
in [Ta2, p. 331]

\begin{lem}
$U_{\mathcal{O}}'/U_{\mathcal{O}}' \cap
E^{-1}\delta^{\vee}\delta(U_{\mathcal{O}})$ is a congruence module
for $E=A_UB_VC^2$.
\end{lem}

$Proof$. Suppose $\phi \in H_{\mathcal{O}}$ acts trivially on
$V_{\mathcal{O}}^{\text{new}}$. We must show that $E\phi$ maps
$U_{\mathcal{O}}'$ into $\delta^{\vee}\delta(U_{\mathcal{O}})$.
Note first that $\phi^{\vee}$ also acts trivially on
$V_{\mathcal{O}}^{\text{new}}$, so it maps $V_{\mathcal{O}}$ into
$V_{\mathcal{O}}^{\text{old}}$. Now let $u =\delta^{\vee}(v) \in
U_{\mathcal{O}}$ for some $v \in V^{\text{old}}$. Note that
$$
A_UC\la v,V_{\mathcal{O}}^{\text{old}} \ra_V \subset A_U\la
v,\delta(U_{\mathcal{O}}) \ra_V \subset A_U\la u,U_{\mathcal{O}}
\ra_U \subset \mathcal{O},
$$
so $A_UC\la \phi v,V_{\mathcal{O}} \ra_V \subset \mathcal{O}
\Rightarrow A_UB_VC(\phi v) \in V_{\mathcal{O}}^{\text{old}}$. We
deduce that
$$
A_UB_VC^2(\phi v) \in \delta(U_{\mathcal{O}}),
$$
and we get the result by applying $\delta^{\vee}$ to this: $E(\phi
u)\in \delta^{\vee}\delta(U_{\mathcal{O}})$. $\qedsymbol$

$\e$

As in [Ta2, p. 331], we have the following useful corollary:

\begin{cor}
Let $\mathcal{O}=\mathcal{O}_L$ be the ring of integers of a
number field $L \subset \C$. Suppose $u \in U_{\mathcal{O}}-\{0\}$
is an eigenvector for $H_{\mathcal{O}}$, with character
$\eta:H_{\mathcal{O}} \rightarrow \mathcal{O}$. Assume:

\begin{itemize}
\item $\mathcal{E} (Lu \cap (U_{\mathcal{O}} + \ker \delta))
\subset \mathcal{O}u$, for some nonzero ideal $\mathcal{E} \subset
\mathcal{O}$, \\
\item $\delta^{\vee}\delta(u) \in m U_{\mathcal{O}}$, for some
nonzero $m \in \mathcal{O}$.
\end{itemize}

Then $\eta$ induces a homomorphism $\T_{\mathcal{O}}^{\text{new}}
\rightarrow \mathcal{O}/\mathcal{O} \cap mE^{-1}\mathcal{E}^{-1}$,
where $E=A_UB_VC^2$.
\end{cor}

We remark that $m=0 \Rightarrow u \in \ker\delta$. If we factor
the fractional ideal $\mathcal{O} \cap mE^{-1}\mathcal{E}^{-1}$
into prime powers and project further, we get the following: For
every (nonzero) prime ideal $\lambda \subset \mathcal{O}$ there is
a homomorphism
$$
\T_{\mathcal{O}}^{\text{new}} \rightarrow \mathcal{O}/\lambda^n
$$
induced by $\eta$, where $n$ is a non-negative integer satisfying
the inequality
$$
n \geq v_{\lambda}(m)-v_{\lambda}(E)-v_{\lambda}(\mathcal{E}).
$$
Here we should think of $v_{\lambda}(m)$ as the main term, and the
other two as controllable error terms. In our applications we want
to show that the right-hand-side is positive.

\section{Compactness at Infinity}

Let $F$ be a totally real number field, and let $\infty$ be the
set of archimedean places. We denote the ring of adeles by
$\A=\A_F=F_{\infty}\times \A^{\infty}$. We consider a connected
reductive $F$-group $G$, and let $A=A_G$ denote the $F$-split
component of its center $Z=Z_G$. Each $F$-rational character $\chi
\in X^*(G)_F$ gives a continuous homomorphism $G(\A)\rightarrow
\R_+^*$ by composing with the idele norm, and we define
$$
G(\A)^1=\{g \in G(\A):|\chi(g)|=1, \forall \chi \in X^*(G)_F\}.
$$
It is known to be unimodular. By the product formula, $G(F)$ is a
discrete subgroup of $G(\A)^1$, and the quotient $G(F) \backslash
G(\A)^1$ has finite volume. Moreover, this quotient is compact if
and only if $G^{\text{ad}}$ is anisotropic. Later, we are
naturally led to studying groups for which
$G_{\infty}^1=G_{\infty}\cap G(\A)^1$ is compact.

\begin{prop}
$G_{\infty}^1$ is compact if and only if one of the following
holds:

\begin{itemize}
\item $G_{\infty}$ is compact, \\
\item $F=\Q$, $\rk_{\Q}Z=\rk_{\R}Z$, and $G_{\infty}^{\text{der}}$
is compact.
\end{itemize}

\end{prop}

$Proof$. Suppose first that $G_{\infty}^1$ is compact. We may
assume that $A \neq 1$ (otherwise $G_{\infty}=G_{\infty}^1$ is
compact). Choosing a basis for $X^*(A)$, we see that (with $r=\dim
A$)
$$
A_{\infty}^1 \simeq \{x \in F_{\infty}^*: \prod_{v\in \infty}
|x_v|_v=1\}^r.
$$
Therefore $\{x \in F_{\infty}^*: \prod_{v|\infty} |x_v|_v=1\}$ is
compact, and we conclude that $F$ has a unique infinite place.
That is, $F=\Q$. If $\rk_{\Q}Z < \rk_{\R}Z$, the $\Q$-anisotropic
component $A'$ is not $\R$-anisotropic. The converse is clear.
$\qedsymbol$

\section{Hecke Algebras and Algebraic Modular Forms}

\subsection{Hecke Algebras}

From now on we fix a totally real number field $F$, and a
connected reductive $F$-group $G$, not a torus, such that
$G_{\infty}^1$ is compact. We consider the locally profinite group
of finite adeles $G(\A^{\infty})$, and choose a Haar measure
$\mu=\otimes \mu_v$ on it once and for all. We consider the vector
space of all locally constant compactly supported $\C$-valued
functions
$$
\mathcal{H}=\mathcal{H}(G(\A^{\infty}))=C_c^{\infty}(G(\A^{\infty}),\C).
$$
This becomes an associative $\C$-algebra, without neutral element,
under $\mu$-convolution. There is a canonical involution on
$\mathcal{H}$ defined by $\phi^{\vee}(x)=\phi(x^{-1})$. This is an
anti-automorphism. If $K \subset G(\A^{\infty})$ is a compact open
subgroup,
$$
e_K =\mu(K)^{-1}\chi_K \in \mathcal{H}
$$
is an idempotent. This is a neutral element in the subalgebra of
$K$-biinvariant compactly supported functions:
$$
\mathcal{H}_K=\mathcal{H}(G(\A^{\infty}),K)=C_c(G(\A^{\infty})//K,\C)=e_K
\star \mathcal{H} \star e_K.
$$
Clearly $\vee$ preserves $\mathcal{H}_K$. In addition, there is a
canonical $\Z$-order $\mathcal{H}_{K,\Z} \subset \mathcal{H}_K$
consisting of all $\mu(K)^{-1}\Z$-valued functions. As a ring,
$\mathcal{H}_{K,\Z}$ is isomorphic to $C_c(G(\A^{\infty})//K,\Z)$
endowed with the $K$-normalized convolution. If $R$ is a
commutative ring, with neutral element, we then define
$$
\mathcal{H}_{K,R}=R \otimes_{\Z}\mathcal{H}_{K,\Z}.
$$
The algebras $\mathcal{H}_K$ are not always commutative. However,
by a result of Bernstein, $\mathcal{H}_K$ is a finite module over
its center $Z(\mathcal{H}_K)$. Now, suppose $J \subset K$ is a
(proper) compact open subgroup. Then obviously $\mathcal{H}_K
\subset \mathcal{H}_J$. However, $\mathcal{H}_K$ is not a subring
since $e_K \neq e_J$. There is a natural retraction $\mathcal{H}_J
\twoheadrightarrow \mathcal{H}_K$ defined by $\phi \mapsto e_K
\star \phi \star e_K$. It does map $e_J \mapsto e_K$, but does not
preserve $\star$ unless we restrict it to the centralizer
$Z_{\mathcal{H}_J}(e_K)$. Clearly,
$Z_{\mathcal{H}_J}(\mathcal{H}_K)$ maps to the center
$Z(\mathcal{H}_K)$. In particular,
$$
\text{$Z(\mathcal{H}_J) \rightarrow Z(\mathcal{H}_K)$, $\y$ $\phi
\mapsto \phi \star e_K = e_K \star \phi$,}
$$
gives a canonical homomorphism of algebras. It maps
$Z(\mathcal{H}_{J,\Z})$ into $Z(\mathcal{H}_{K,\Z})$.

\subsection{Algebraic Modular Forms}

Note that $G(F) \subset G(\A^{\infty})$ is a discrete subgroup. We
consider the Hilbert space of $L^2$-functions on the quotient,
$L^2(G(F) \backslash G(\A^{\infty}))$. There is a unitary
representation $r$ of $G(\A^{\infty})$ on this space given by
right translations. We consider the smooth vectors,
$$
\mathcal{A}=L^2(G(F) \backslash
G(\A^{\infty}))^{\infty}=C^{\infty}(G(F) \backslash
G(\A^{\infty}),\C),
$$
consisting of locally constant functions. This is an admissible
representation:
$$
\text{$\mathcal{A}=\bigcup \mathcal{A}_K$, $\y$ where $\y$
$\mathcal{A}_K \simeq C(G(F) \backslash G(\A^{\infty})/K,\C)$,}
$$
and $K$ runs over all compact open subgroups of $G(\A^{\infty})$.
Therefore, the Hecke algebra $\mathcal{H}$ acts on $\mathcal{A}$
in the usual way. We have the following compatibility between this
action and the inner product:
$$
(r(\bar{\phi})f,g)=(f,r(\phi^{\vee})g).
$$
For a compact open subgroup $K \subset G(\A^{\infty})$, the space
of $K$-invariants
$$
\mathcal{A}_K \simeq C(G(F) \backslash
G(\A^{\infty})/K,\C)=r(e_K)\mathcal{A}
$$
is finite-dimensional. Indeed the double coset space $X_K=G(F)
\backslash G(\A^{\infty})/K$ is finite. Functions in
$\mathcal{A}_K$ are examples of algebraic modular forms. Clearly,
$\mathcal{H}_K$ acts on $\mathcal{A}_K$, and the order
$\mathcal{H}_{K,\Z}$ preserves the lattice of $\Z$-valued
functions:
$$
\mathcal{A}_{K,\Z}=C(G(F) \backslash G(\A^{\infty})/K,\Z) \subset
\mathcal{A}_K.
$$
For a commutative ring $R$ we let $\mathcal{A}_{K,R}=R
\otimes_{\Z}\mathcal{A}_{K,\Z}$. The $R$-algebra
$\mathcal{H}_{K,R}$ acts on this module, and we let $\T_{K,R}$
denote the image of the center $Z(\mathcal{H}_{K,R})$ in $\End_R
\mathcal{A}_{K,R}$. Hence $\T_{K,R}$ is a commutative $R$-algebra.
Now, suppose $J \subset K$ is a (proper) compact open subgroup.
Then $\mathcal{A}_K \subset \mathcal{A}_J$, and the canonical
homomorphism $Z(\mathcal{H}_{J,R}) \rightarrow
Z(\mathcal{H}_{K,R})$ descends to the restriction map
$\T_{J,R}\rightarrow \T_{K,R}$.

\subsection{Pairings}

We define a pairing on $\mathcal{A}_K$ as follows. Here $(-,-)$
denotes the Petersson inner product.

\begin{df}
For $f, g \in \mathcal{A}_K$, we define a symmetric bilinear form
by
$$
\la f,g \ra_K=\mu(K)^{-1}(f,\bar{g})=\sum_{x \in X_K}
f(x)g(x)|G(F)\cap {^xK}|^{-1},
$$
where we use the notation ${^xK}=xKx^{-1}$.
\end{df}

The factors $|G(F)\cap {^xK}|^{-1}$ are missing in [Ta1] and
[Ta2]. If $K$ is sufficiently small, for example if $K=\prod_{v <
\infty}K_v$ and some $K_v$ is torsion-free (this is the case if
$K_v$ is a sufficiently deep principal congruence subgroup), then
indeed $G(F)\cap {^xK}=1$. For $\phi \in \mathcal{H}_K$ and $f, g
\in \mathcal{A}_K$ we have the compatibility relation
$$
\la r(\phi)f,g \ra_K=\la f, r(\phi^{\vee})g \ra_K.
$$
Next we have to show the quotient
$\mathcal{A}_{K,\OO}/\mathcal{A}_{K,\OO}^{\vee}$ is torsion and
find a good annihilator $A_K$. The fact that it is torsion is
immediate: It is killed by the positive integer
$$
\prod_{x \in X_K}|G(F)\cap {^xK}|.
$$
This is $1$ if $K$ is sufficiently small in the sense above.

\begin{lem}
Let $K=\prod_{v<\infty}K_v \subset G(\A^{\infty})$ be a
decomposable compact open subgroup, and let $\ell$ be a prime
number. Suppose $\ell \nmid |K_v|$ for some $v < \infty$. Then
there exists a positive integer $A_K$, not divisible by $\ell$,
such that
$$
A_K \la \mathcal{A}_{K,\OO},\mathcal{A}_{K,\OO} \ra_K \subset \OO.
$$
\end{lem}

$Proof$. Choose some torsion-free subgroup $\tilde{K}_v \subset
K_v$ and let $\tilde{K}=\tilde{K}_vK^v$. Then
$$
\la \mathcal{A}_{\tilde{K},\OO},\mathcal{A}_{\tilde{K},\OO}
\ra_{\tilde{K}} \subset \OO
$$
as we have observed above. Therefore, for $f, g \in
\mathcal{A}_{K,\OO} \subset \mathcal{A}_{\tilde{K},\OO}$, we have
$$
[K_v:\tilde{K}_v]\la f,g \ra_K=\la f,g \ra_{\tilde{K}} \in \OO.
$$
We then take $A_K=[K_v:\tilde{K}_v]$. This is not divisible by
$\ell$. $\qedsymbol$

$\e$

Note that $\ell \nmid |K_v|$ if $K_v$ is torsion-free and $v \nmid
\ell$. For large $\ell$ this is automatic:

\begin{lem}
Suppose there exists an $F$-embedding $G \hookrightarrow \GL(n)$.
Let $K=\prod_{v<\infty}K_v$ be arbitrary and let $\ell >
[F:\Q]n+1$ be a prime number. Then $\ell \nmid |K_v|$ holds for
infinitely many places $v$.
\end{lem}

$Proof$. $K_v$ embeds into a conjugate of
$\GL(n,\mathcal{O}_{v})$. Therefore $|K_v|$ divides
$$
|\GL(n,\mathcal{O}_{v})|=p^{\infty}\prod_{i=1}^n(q^i-1).
$$
Assume $\ell$ divides $|K_v|$ for almost all $v$. Then $p$ has
order at most $[F:\Q]n$ in $(\Z/\ell)^*$ for almost all primes
$p$. Now, $(\Z/\ell)^*$ is cyclic of order $\ell-1$, so by
Dirichlet's theorem on primes in arithmetic progressions we
conclude that $\ell \leq [F:\Q]n+1$. $\qedsymbol$

\section{Parahoric Level Structure and the Concrete Setup}

\subsection{Parahoric Subgroups}

From now on we assume for simplicity that $G^{\text{der}}$ is
simple (that is, it has no nontrivial connected normal subgroups).
Moreover, we fix a compact open subgroup
$$
K=\prod_{v < \infty}K_v \subset G(\A^{\infty}).
$$
It is known that $K_v \subset G_v$ is a hyperspecial maximal
compact subgroup for almost all places $v$, that is,
$K_v=\underline{G}(\mathcal{O}_{F_v})$ for a smooth affine group
scheme $\underline{G}$ of finite type over $\mathcal{O}_{F_v}$
with generic fiber $G$. Such exist precisely when $G_v$ is
unramified. Let us look at a fixed finite place $w$ of $F$ where
$K_w$ is hyperspecial. Then write $K=K_wK^w$, where
$$
K^w=\prod_{v \neq w}K_v \subset G(\A^{\infty,w}).
$$
Let $\mathcal{B}_w$ denote the reduced Bruhat-Tits building of
$G_w$ (that is, the building of $G_w^{\text{ad}}$). We have
assumed $G^{\text{der}}$ is simple, so $\mathcal{B}_w$ is a
simplicial complex. Let $x \in \mathcal{B}_w$ be the vertex fixed
by $K_w$, and let $(x,x')$ be an edge in the building. Then
consider the maximal compact subgroup $K_w' \subset G_w$ fixing
the vertex $x'$, and the parahoric subgroup $J_w=K_w \cap K_w'$
associated with the edge $(x,x')$. Let $K'=K_w'K^w$ and $J=J_wK^w$
be the corresponding subgroups of $G(\A^{\infty})$.

\begin{lem}
$\la K_w,K_w' \ra=G_w^0:=\{g \in G_w: |\chi(g)|=1, \forall \chi
\in X^*(G)_{F_w}\}$.
\end{lem}

$Proof$. This follows from Bruhat-Tits theory. $\qedsymbol$

$\e$

Note that $G_w^{\text{der}} \subset G_w^0 \subset G_w^1=G_w \cap
G(\A)^1$.

\subsection{The Concrete Setup}

Now we want to apply our general results in the following setup:
Let $L \subset \C$ be a number field, and let
$\mathcal{O}=\mathcal{O}_L$ be its ring of integers. We let
$H=Z(\mathcal{H}_J)$. This is a commutative $\C$-algebra, and it
comes with the involution defined by
$\phi^{\vee}(x)=\phi(x^{-1})$. $V=\mathcal{A}_J$ is a
finite-dimensional $\C$-space on which $Z(\mathcal{H}_J)$ acts.
The order $Z(\mathcal{H}_{J,\mathcal{O}})$ preserves the lattice
$V_{\mathcal{O}}=\mathcal{A}_{J,\mathcal{O}}$. The space $V$ comes
with the bilinear form $\la -,- \ra_J$. The compatibility
conditions between these data are satisfied. Let $U=\mathcal{A}_K
\oplus \mathcal{A}_{K'}$. Then $Z(\mathcal{H}_J)$ acts on this
space via the natural maps to $Z(\mathcal{H}_K)$ and
$Z(\mathcal{H}_{K'})$. The lattice $U_{\mathcal{O}}=
\mathcal{A}_{K,\mathcal{O}} \oplus \mathcal{A}_{K',\mathcal{O}}$
is preserved by $Z(\mathcal{H}_{J,\mathcal{O}})$. The bilinear
form on $U$ is given by the sum $\la -,- \ra_K \oplus \la -,-
\ra_{K'}$. The degeneracy map $\delta$ is given by
$$
\delta: \mathcal{A}_K \oplus \mathcal{A}_{K'}
\overset{\text{sum}}{\rightarrow} \mathcal{A}_J,
$$
which is clearly $Z(\mathcal{H}_J)$-linear. Obviously,
$\ker(\delta)$ consists of all pairs $(f,-f)$, where
$$
f \in \mathcal{A}_K \cap
\mathcal{A}_{K'}=\{\text{$G_w^0K^w$-invariant functions $f \in
\mathcal{A}$}\}.
$$
The decompositions $U=\ker(\delta)\oplus\ker(\delta)^{\bot}$ and
$V=\im(\delta)\oplus\im(\delta)^{\bot}$ are immediate because of
the relation between the pairings and the inner product.

\subsection{The Combinatorial Ihara Lemma}

The proof of the following lemma is a straightforward
generalization of [Ta1, p. 274]:

\begin{lem}
$\mathcal{A}_{J,\mathcal{O}} \cap \delta(\mathcal{A}_K \oplus
\mathcal{A}_{K'})=\delta(\mathcal{A}_{K,\mathcal{O}} \oplus
\mathcal{A}_{K',\mathcal{O}})$.
\end{lem}

$Proof$. Let us first set up some machinery for the proof. We
define an equivalence relation on $X_J$ by saying that $x,y \in
X_J$ are equivalent ($x \sim y$) if and only if
$$
\text{$\exists$ chain $x=x_0,\ldots,x_d=y$ such that $\forall i$:
$\pi(x_i)=\pi(x_{i+1})$ or $\pi'(x_i)=\pi'(x_{i+1})$.}
$$
This gives a partition of $X_J$ into equivalence classes $X_J^j$.
For each $j$, we fix a representative $y^j \in X_J^j$.
Correspondingly, we have a radius function $d:X_J \rightarrow
\Z_{\geq 0}$ defined as follows: Given $x \in X_J$, there is a
unique $j$ such that $x \sim y^j$. Then $d(x)$ is the minimal
length of a chain connecting $x$ to $y^j$. Now, suppose
$g=\delta(f,f') \in \mathcal{A}_{J,\mathcal{O}}$ for some $f \in
\mathcal{A}_{K}$ and $f' \in \mathcal{A}_{K'}$. We want to show $g
\in \delta(\mathcal{A}_{K,\mathcal{O}} \oplus
\mathcal{A}_{K',\mathcal{O}})$.
$$
\text{\underline{Claim} - We may assume that $f(\pi(y^j))=0$ for
all $j$.}
$$
To see this, note that $X_K=\sqcup \pi(X_J^j)$ and $X_{K'}=\sqcup
\pi'(X_J^j)$. We then define $\tilde{f} \in \mathcal{A}_{K}$ such
that $\tilde{f}|\pi(X_J^j)\equiv f(\pi(y^j))$, and $\tilde{f}' \in
\mathcal{A}_{K'}$ such that $\tilde{f}'|\pi'(X_J^j)\equiv
f(\pi(y^j))$. Then
$$
g=\delta(f-\tilde{f},f'+\tilde{f}')
$$
and $(f-\tilde{f})(\pi(y^j))=0$ for all $j$. This proves the
claim, so from now on assume that $f(\pi(y^j))=0$ for all $j$. By
induction on $m \geq 0$, we now prove the following:
$$
\text{$\forall x \in X_J$ with $d(x)=m$: $f(\pi(x)) \in
\mathcal{O}$ and $f'(\pi'(x)) \in \mathcal{O}$.}
$$
This is sufficient, for then $f \in \mathcal{A}_{K,\mathcal{O}}$
and $f' \in \mathcal{A}_{K',\mathcal{O}}$. Note that $f(\pi(x))
\in \mathcal{O}$ if and only if $f'(\pi'(x)) \in \mathcal{O}$. The
start $m=0$ is essentially just our assumption, so assume the
statement is true for $m-1 \geq 0$ and consider $x \in X_J$ with
$d(x)=m$. Let
$$
x=x_0,x_1,\ldots,x_m=y^j
$$
be a chain of minimal length. Then $x'=x_1 \in X_J$ has
$d(x')=m-1$, so by induction $f(\pi(x')) \in \mathcal{O}$ and
$f'(\pi'(x')) \in \mathcal{O}$. However, $\pi(x)=\pi(x')$ or
$\pi'(x)=\pi'(x')$. In either case we get the statement for $x$.
$\qedsymbol$

\section{Applying the Abstract Theory}

\subsection{Computing $\delta^{\vee}\delta$}

To apply the abstract theory it is necessary to compute
$\delta^{\vee}\delta$ explicitly.

\begin{lem} The endomorphism $\delta^{\vee}\delta$ is given by the
$2\times 2$ matrix
$$
\delta^{\vee}\delta=\begin{pmatrix}[K:J] & [K:J]e_K \\
[K':J]e_{K'} & [K':J] \end{pmatrix}.
$$
\end{lem}

$Proof$. $\delta^{\vee}\delta$ is an endomorphism of
$\mathcal{A}_K \oplus \mathcal{A}_{K'}$, and we write is as
$$
\delta^{\vee}\delta=\begin{pmatrix}a & b \\ c & d\end{pmatrix},
$$
where $b: \mathcal{A}_{K'} \rightarrow \mathcal{A}_{K}$ and so on.
Using the definition it is not hard to see that
$$
\la af,g \ra_K=\la f,g\ra_J=[K:J]\la f,g\ra_K
$$
for all $f,g \in \mathcal{A}_{K}$. In particular, $a=[K:J]$. In
the same way one computes $b$, $c$ and $d$. $\qedsymbol$

\subsection{The Main Lemma}

In our situation, Corollary 1 gives the following crucial lemma.

\begin{lem}
Let $f \in \mathcal{A}_{K,\mathcal{O}}$ be an eigenform for
$Z(\mathcal{H}_{K,\mathcal{O}})$ with character
$\eta_f:\T_{K,\mathcal{O}} \rightarrow \mathcal{O}$. Assume $f$ is
not $G_w^0$-invariant modulo $\lambda|\ell$, where $\ell$
satisfies the following:
$$
\text{There exists at least two places $v$ such that $\ell \nmid
|K_v|$.}
$$
Then the reduction of $\eta_f$ modulo $\lambda^n$ factors through
$\T_{J,\mathcal{O}}^{\text{new}}$ when
$$
v_{\lambda}(\eta_f(e_{K,K'})-[K:J][K':J])-v_{\lambda}([K':J]) \geq
n,
$$
where we introduce the notation $e_{K,K'}=[K:J][K':J](e_K \star
e_{K'}\star e_K) \in Z(\mathcal{H}_{K,\Z})$.
\end{lem}

$Proof$. To produce an eigenvector in $U_{\mathcal{O}}=
\mathcal{A}_{K,\mathcal{O}} \oplus \mathcal{A}_{K',\mathcal{O}}$,
we take
$$
\vec{f}=[K':J](f,-r(e_{K'})f) \in \mathcal{A}_{K,\mathcal{O}}
\oplus \mathcal{A}_{K',\mathcal{O}}.
$$
The factor $[K':J]$ is included since $r(e_{K'})f$ does not
necessarily take values in $\mathcal{O}$. Clearly, $\vec{f}$ is an
eigenvector for $Z(\mathcal{H}_{J,\mathcal{O}})$, and its
character is the composite
$$
\eta_{\vec{f}}: Z(\mathcal{H}_{J,\mathcal{O}}) \rightarrow
Z(\mathcal{H}_{K,\mathcal{O}}) \overset{\eta_f}{\rightarrow}
\mathcal{O}.
$$
Using the explicit formula for $\delta^{\vee}\delta$ in lemma 6
above, it follows that
$$
\delta^{\vee}\delta(\vec{f})=(\eta_f(e_{K,K'})-[K:J][K':J])(-f,0).
$$
Now, since $(-f,0) \in U_{\mathcal{O}}$, in Corollary 1 we can
take
$$
m=\eta_f(e_{K,K'})-[K:J][K':J] \in \mathcal{O}
$$
as long as this is nonzero. However, note that $\vec{f}$ must
belong to the kernel of $\delta$ if $m=0$. Hence $f$ must be
invariant under the group $G_w^0$ (say, on the right). Now, let
$$
\mathcal{F}=\{x \in L: xf \in
\mathcal{A}_{K,\mathcal{O}}+\mathcal{A}_K\cap\mathcal{A}_{K'}\}.
$$
This is clearly an $\mathcal{O}$-submodule of $L$ containing
$\mathcal{O}$. Obviously, $\mathcal{F}=L$ if $f \in
\mathcal{A}_K\cap\mathcal{A}_{K'}$. However, $f$ is not
$G_w^0$-invariant, so $\mathcal{F}$ is a fractional ideal. To see
this note that
$$
A_K \la f,g \ra_K \mathcal{F} \subset \mathcal{O}
$$
for every $g \in \mathcal{A}_{K,\mathcal{O}} \cap
(\mathcal{A}_K\cap\mathcal{A}_{K'})^{\bot}$. These $g$ span
$(\mathcal{A}_K\cap\mathcal{A}_{K'})^{\bot}$ so $f$ must belong to
$\mathcal{A}_K\cap\mathcal{A}_{K'}$ if $\la f,g \ra_K=0$ for all
such $g$. Now, the nonzero ideal
$\tilde{\mathcal{E}}=\mathcal{F}^{-1}$ satisfies:
$$
\tilde{\mathcal{E}}(Lf \cap
(\mathcal{A}_{K,\mathcal{O}}+\mathcal{A}_K\cap\mathcal{A}_{K'}))
\subset \mathcal{O}f.
$$
Therefore, $\mathcal{E}=[K':J]\tilde{\mathcal{E}}$ satisfies the
primitivity condition in corollary 1:
$$
\mathcal{E}(L\vec{f} \cap (\mathcal{A}_{K,\mathcal{O}} \oplus
\mathcal{A}_{K',\mathcal{O}}+\ker \delta)) \subset
\mathcal{O}\vec{f}.
$$
Suppose $\lambda \subset \mathcal{O}$ is a maximal ideal such that
$v_{\lambda}(\tilde{\mathcal{E}})\neq 0$. Then $\lambda^{-1}
\subset \mathcal{F}$. It follows that $f \in
\lambda(\mathcal{A}_{K,\mathcal{O}}+\mathcal{A}_K\cap\mathcal{A}_{K'})$,
and hence the reduction $\bar{f} \in \mathcal{A}_{K,\F_{\lambda}}$
is $G_w^0$-invariant. Since $\ell \nmid |K_v|$ holds for at least
one $v \neq w$, by assumption, we can find $A_K$ and $A_{K'}$
indivisible by $\ell$ according to Lemma 2. Note also that we can
take $C=1$ by Lemma 5. $\qedsymbol$

\section{Semisimplicity}

\subsection{Automorphic Representations and the Decomposition of $\mathcal{A}_K$}

Henceforth assume $G^{\text{der}}$ is simple and simply connected.
There is an admissible representation of $G(\A^{\infty})$ on the
space
$$
\mathcal{A}=C^{\infty}(G(F)\backslash G(\A^{\infty}),\C)
$$
given by right translations. For a compact open subgroup $K
\subset G(\A^{\infty})$, we look at the $\mathcal{H}_K$-module of
$K$-invariants $\mathcal{A}_K$. Recall that an automorphic
representation of $G(\A)$ is an irreducible representation $\pi$
of $G(\A)$ (on some Hilbert space) such that
$$
\Hom_{G(\A)^1}(\pi,L^2(G(F)\backslash G(\A)^1))\neq 0.
$$
We let $m(\pi)$ denote the dimension of this space. We then have
an isomorphism,
$$
\mathcal{A}_K \simeq {\bigoplus}_{\pi \in \Pi_{\text{unit}}(G(\A))
: \pi_{\infty}={\bf{1}}} m(\pi)\pi^K,
$$
On the right we have a finite direct sum over the automorphic
representations $\pi$ of $G(\A)$ such that $\pi_{\infty}={\bf{1}}$
and $\pi^K \neq 0$. These $\pi$ are automatically unitary.

\subsection{Semisimplicity in Characteristic Zero}

It is known that each $\pi^K$ is a simple module over
$\mathcal{H}_K$, and hence $\mathcal{A}_K$ is semisimple.
Moreover, by Schur's lemma, the center $Z(\mathcal{H}_K)$ acts on
$\pi^K$ by a $\C$-algebra homomorphism
$\eta_{\pi^K}:Z(\mathcal{H}_K) \rightarrow \C$. For a character
$\eta:Z(\mathcal{H}_K) \rightarrow \C$, we denote by
$\mathcal{A}_K(\eta)$ the $\eta$-isotypic component. That is, the
eigenspace
$$
\mathcal{A}_K(\eta)=\{f \in \mathcal{A}_K: r(\phi)f=\eta(\phi)f,
\forall \phi \in Z(\mathcal{H}_K)\}.
$$
Then there is a direct sum decomposition
$\mathcal{A}_K=\bigoplus_{\eta}\mathcal{A}_K(\eta)$. Clearly,
$\mathcal{A}_K(\eta) \neq 0$ if and only if $\eta=\eta_{\pi^K}$
for some $\pi$. The image $\T_K \subset \End \mathcal{A}_K$ of the
center $Z(\mathcal{H}_K)$ is a commutative semisimple
$\C$-algebra, that is, a direct product of copies of $\C$.

\begin{lem}
The eigenspace $\mathcal{A}_K(\eta)$ is nonzero if and only if
$\eta$ factors through $\T_K$.
\end{lem}

$Proof$. Obviously, $\eta$ factors if $\mathcal{A}_K(\eta)\neq 0$.
Conversely, suppose $\eta$ factors and look at its kernel
$\m=\ker(\eta) \subset \T_K$. This is a maximal ideal. Since
$\T_K$ acts faithfully on $\mathcal{A}_K$, which is
finite-dimensional, $\m$ belongs to the support of
$\mathcal{A}_K$. By the theory of associated primes, $\m$ contains
a prime ideal of the form $\Ann_{\T_K}(f)$ with $f \in
\mathcal{A}_K$. All primes are maximal in $\T_K$, so in fact
$\m=\Ann_{\T_K}(f)$. Clearly $\m$ contains $T-\eta(T)$ for every
$T \in \T_K$, so $f \in \mathcal{A}_K(\eta)$, and $f$ must be
nonzero as $\m \neq \T_K$. $\qedsymbol$

$\e$

Now, consider the $\mathcal{H}_{K,\Q}$-module
$\mathcal{A}_{K,\Q}$, and the image $\T_{K,\Q}$ of the center
$Z(\mathcal{H}_{K,\Q})$ in the endomorphism algebra
$\End_{\Q}\mathcal{A}_{K,\Q}$. $\T_{K,\Q}$ can be viewed as a
subring of $\T_K \simeq \C \otimes_{\Q} \T_{K,\Q}$. We deduce that
$\T_{K,\Q}$ is a reduced commutative finite-dimensional
$\Q$-algebra, that is, a product of number fields by Nakayama's
lemma:
$$
\T_{K,\Q} \simeq L_1 \times \cdots \times L_t.
$$
Visibly, $\T_{K,\Q}$ is a semisimple $\Q$-algebra. The $L_i$
occurring in $\T_{K,\Q}$ are totally real or CM.

\subsection{Semisimplicity in Positive Characteristic}

Now let $R$ be a field of characteristic $p>0$. We are interested
in when $\mathcal{A}_{K,R}$ is a semisimple module over
$Z(\mathcal{H}_{K,R})$. As we have seen, this means that
$\T_{K,R}$ is a semisimple $R$-algebra. We have $\T_{K,R}\simeq R
\otimes_{\F_p}\T_{K,\F_p}$, so equivalently, when is $\T_{K,\F_p}$
semisimple? There is always a surjective homomorphism
$$
\xi:\F_p \otimes_{\Z} \T_{K,\Z} \twoheadrightarrow \T_{K,\F_p}.
$$
Indeed the image of $\F_p \otimes_{\Z} \T_{K,\Z}$ in
$\End_{\F_p}\mathcal{A}_{K,\F_p}$ equals the image of $\F_p
\otimes_{\Z} Z(\mathcal{H}_{K,\Z})$. However, the natural map from
this last algebra to $Z(\mathcal{H}_{K,\F_p})$ is surjective. Let
$$
\tilde{\T}_{K,\Z}=\{T \in \T_{K,\Q}: T(\mathcal{A}_{K,\Z})\subset
\mathcal{A}_{K,\Z}\}.
$$
This is a finite free $\Z$-module containing $\T_{K,\Z}$ as a
subgroup of finite index.

\begin{lem}
The kernel $\ker(\xi)$ is nilpotent. It is trivial if and only if
$p \nmid [\tilde{\T}_{K,\Z}:\T_{K,\Z}]$.
\end{lem}

$Proof$. It is enough to show that every element in $\ker(\xi)$ is
nilpotent. Under the identification $\F_p \otimes_{\Z} \T_{K,\Z}
\simeq \T_{K,\Z}/p\T_{K,\Z}$, the kernel $\ker(\xi)$ corresponds
to the ideal
$$
\T_{K,\Z} \cap p\tilde{\T}_{K,\Z}/p\T_{K,\Z}.
$$
Let $T \in \T_{K,\Z} \cap p\tilde{\T}_{K,\Z}$. Obviously,
$\tilde{\T}_{K,\Z}$ is integral over $\Z$, so there is an equation
$$
(p^{-1}T)^n+a_{n-1}(p^{-1}T)^{n-1}+\cdots+a_1(p^{-1}T)+a_0=0
$$
for certain $a_i \in \Z$. Multiplying by $p^n$ we see that $T^n
\in p\T_{K,\Z}$. For the last assertion, note that $\ker(\xi)=0$
if and only if $\F_p \otimes_{\Z} \T_{K,\Z} \rightarrow \F_p
\otimes_{\Z} \tilde{\T}_{K,\Z}$ is injective. $\qedsymbol$

$\e$

In particular, $\ker(\xi)$ is contained in the Jacobson radical.
We let $\bar{\T}_{K,\Z}$ denote the integral closure of $\Z$ in
$\T_{K,\Q}$. It contains $\tilde{\T}_{K,\Z}$ as a subgroup of
finite index.

\begin{lem}
$p \nmid \Delta_K:=
[\bar{\T}_{K,\Z}:\tilde{\T}_{K,\Z}]\cdot\prod_{i}\Delta_{L_i/\Q}
\Rightarrow \text{$\T_{K,\F_p}$ is semisimple}$.
\end{lem}

$Proof$. Note first that $\F_p \otimes_{\Z}\tilde{\T}_{K,\Z}
\simeq \F_p \otimes_{\Z}\bar{\T}_{K,\Z}$ since $p \nmid
[\bar{\T}_{K,\Z}:\tilde{\T}_{K,\Z}]$. Now,
$$
\F_p \otimes_{\Z}\bar{\T}_{K,\Z} \simeq \prod_i
\mathcal{O}_{L_i}/p\mathcal{O}_{L_i} \simeq \prod_i
\prod_{\p|p}\mathcal{O}_{L_i}/\p,
$$
since $p$ is unramified in every $L_i$ occurring in $\T_{K,\Q}$.
There is an embedding,
$$
\T_{K,\F_p} \simeq \T_{K,\Z}/\T_{K,\Z} \cap p\tilde{\T}_{K,\Z}
\hookrightarrow \tilde{\T}_{K,\Z}/p\tilde{\T}_{K,\Z} \simeq \F_p
\otimes_{\Z}\tilde{\T}_{K,\Z},
$$
and it follows that $\T_{K,\F_p}$ is semisimple. $\qedsymbol$

$\e$

The converse holds at least for $p \nmid
[\tilde{\T}_{K,\Z}:\T_{K,\Z}]$ (that is, when $\xi$ is injective).

\subsection{The Simple Modules}

Let $R$ be a perfect field of characteristic $p \geq 0$. Up to
isomorphism, the simple $Z(\mathcal{H}_{K,R})$-modules are given
by an extension $R'/R$ with an action given by a surjective
$R$-algebra homomorphism $\eta:Z(\mathcal{H}_{K,R})
\twoheadrightarrow R'$. If $\eta$ is a submodule of
$\mathcal{A}_{K,R}$, the extension $R'/R$ is finite and $\eta$
factors through $\T_{K,R}$. If $p \nmid \Delta_K$, there exists a
finite extension $L/R$ such that we have a direct sum
decomposition
$$
\mathcal{A}_{K,L}=\bigoplus_{\eta}\mathcal{A}_{K,L}(\eta),
$$
This is still true when $p|\Delta_K$, provided
$\mathcal{A}_{K,L}(\eta)$ denotes the generalized eigenspace:
$$
\mathcal{A}_{K,L}(\eta)=\{\text{$f \in \mathcal{A}_{K,L}$:
$\forall \phi \in  Z(\mathcal{H}_{K,L})$,
$(r(\phi)-\eta(\phi))^nf=0$ for some $n \geq 1$}\}.
$$
Observe the following:

\begin{lem} Let $R$ be a field, and choose a finite extension
$L/R$ as above. Then let $L'/L$ be an arbitrary extension. Suppose
$\eta': Z(\mathcal{H}_{K,L'}) \rightarrow L'$ occurs in
$\mathcal{A}_{K,L'}$. Then $\eta'=1\otimes \eta$ for some
character $\eta: Z(\mathcal{H}_{K,L}) \rightarrow L$ occurring in
$\mathcal{A}_{K,L}$. Moreover,
$$
\mathcal{A}_{K,L'}(1\otimes \eta) \simeq L' \otimes_L
\mathcal{A}_{K,L}(\eta),
$$
so $\eta$ and $\eta'=1 \otimes \eta$ occur with the same
multiplicity.
\end{lem}

$Proof$. Both $\mathcal{A}_{K,L}$ and $\mathcal{A}_{K,L'} \simeq
L' \otimes_L \mathcal{A}_{K,L}$ have decompositions into direct
sums of generalized eigenspaces. Under this isomorphism, $L'
\otimes_L \mathcal{A}_{K,L}(\eta) \hookrightarrow
\mathcal{A}_{K,L'}(1\otimes \eta)$. Therefore, every $\eta'$
occurring in $\mathcal{A}_{K,L'}$ must come from an $\eta$, and
the above injection must be an isomorphism. $\qedsymbol$

$\e$

Let us apply these results to $R=\Q$. We conclude that there
exists a number field $L/\Q$ such that $\mathcal{A}_{K,L}$ is a
direct sum of eigenspaces for characters $Z(\mathcal{H}_{K,L})
\rightarrow L$. Furthermore, if $\eta: Z(\mathcal{H}_K)
\rightarrow \C$ is a character such that $\mathcal{A}_K(\eta)\neq
0$, then $\eta$ restricts to a $\Q$-algebra homomorphism
$Z(\mathcal{H}_{K,\Q}) \rightarrow L$ occurring in
$\mathcal{A}_{K,L}$. In addition, since $Z(\mathcal{H}_{K,\Z})$
preserves $\mathcal{A}_{K,\mathcal{O}_L}$, $\eta$ even restricts
to a ring homomorphism $Z(\mathcal{H}_{K,\Z}) \rightarrow
\mathcal{O}_L$ occurring in $\mathcal{A}_{K,\mathcal{O}_L}$.

\section{End of the Proof}

\subsection{Invariance Modulo $\lambda$}

The following is a more refined version of the notion abelian
modulo $\lambda$.

\begin{df}
Let $\pi$ be an automorphic representation of $G(\A)$ such that
$\pi_{\infty}={\bf{1}}$. We say that $\pi$ is abelian modulo
$\lambda$ relative to $K$ if $\pi^K \neq 0$ and there exists an
automorphic character $\chi$ of $G(\A)$, trivial on $G_{\infty}K$,
such that $ \text{$\eta_{\pi^K}(\phi) \equiv \eta_{\chi}(\phi)$
(mod $\lambda$)}$, $\forall \phi \in Z(\mathcal{H}_{K,\Z})$.
\end{df}

If this holds, we can find eigenforms in $\pi^K$ to which our main
lemma applies:

\begin{lem}
Let $\pi$ be an automorphic representation of $G(\A)$ such that
$\pi_{\infty}={\bf{1}}$. If $\pi$ is non-abelian modulo $\lambda$
relative to $K$, then the eigenspace
$\mathcal{A}_{K,\F}^{\circ}(\bar{\eta})$ contains no nonzero
$G_w^{\text{der}}$-invariant functions, where $w$ is a place such
that $K_w$ is hyperspecial.
\end{lem}

$Proof$. Choose a number field $L/\Q$ such that
$\mathcal{A}_{K,L}$ is a direct sum of eigenspaces and let
$\mathcal{O}=\mathcal{O}_L$. Denote by
$$
\eta=\eta_{\pi^K}: Z(\mathcal{H}_{K,\Z}) \rightarrow \mathcal{O}
$$
the character giving the action on $\pi^K$. As we have observed
above it occurs in $\mathcal{A}_{K,\mathcal{O}}$, that is, there
exists an eigenform $0 \neq f \in \mathcal{A}_{K,\mathcal{O}}$
with $\eta_f=\eta$. We consider a finite place $\lambda|\ell$ of
$\bar{\Q}$, and a finite place $w$ of $F$ such that $G_w$ is
unramified. Let $\bar{f}=1 \otimes f \in \mathcal{A}_{K,\F}$ be
the reduction modulo $\lambda$, where $\F=\mathcal{O}/\lambda \cap
\mathcal{O}$ is a finite extension of $\F_{\ell}$. By scaling $f$,
we can assume that $\bar{f} \neq 0$. Let us assume $\bar{f}$ is
$G_w^{\text{der}}$-invariant. Now, $G^{\text{der}}$ is simple,
simply connected and $G_w^{\text{der}}$ is noncompact. By the
strong approximation theorem, $\bar{f}$ is in fact
$G^{\text{der}}(\A^{\infty})$-invariant. There is a short exact
sequence
$$
1 \rightarrow G^{\text{der}}(\A^{\infty}) \rightarrow
G(\A^{\infty})
\overset{\nu}{\rightarrow}G^{\text{ab}}(\A^{\infty}) \rightarrow
1.
$$
It follows that $\bar{f}$ lives on $G^{\text{ab}}(\A^{\infty})$.
More precisely, there exists a unique function
$\tilde{f}:G^{\text{ab}}(\A^{\infty}) \rightarrow \F$ such that
$\bar{f}=\tilde{f} \circ \nu$. It fits into the diagram
\begin{align*}
\xymatrix{ X_K=G(F) \backslash G(\A^{\infty})/K \ar[r]^-{\bar{f}}
\ar@{->>}[d]_{\nu} &
\F \\
Y_K=\nu(G(F))\backslash G^{\text{ab}}(\A^{\infty})/\nu(K)
\ar[ur]_-{\tilde{f}} }
\end{align*}

If $R$ is a ring we denote by $\mathcal{A}_{K,R}^{\text{ab}}$ the
module of $R$-valued functions on $Y_K$. Pulling back via $\nu$,
identifies $\mathcal{A}_{K,R}^{\text{ab}}$ with an
$\mathcal{H}_{K,R}$-submodule of $\mathcal{A}_{K,R}$. Then $0 \neq
\tilde{f} \in \mathcal{A}_{K,\F}^{\text{ab},
  \circ}(\bar{\eta})$.
By the Deligne-Serre lifting lemma (that is, Lemme 6.11 in their
paper [DS, p. 522]) we can lift $\bar{\eta}$ to characteristic
zero: There exists an eigenform $0 \neq f' \in
\mathcal{A}_{K,L_{\lambda}}^{\text{ab}}$ such that its character
$\eta': Z(\mathcal{H}_{K,\Z}) \rightarrow \mathcal{O}_{\lambda}$
reduces to $\bar{\eta}$ modulo $\lambda \cap \mathcal{O}$. From
the results of the previous section we see that in fact $\eta'$
maps into $\mathcal{O}$, and it occurs in
$\mathcal{A}_{K,L}^{\text{ab}}$ (and therefore in
$\mathcal{A}_{K}^{\text{ab}}$). However,
$\mathcal{A}_{K}^{\text{ab}}$ is just the space of $\C$-valued
functions on the finite abelian group $Y_K$, so the characters
form a basis. We conclude that there exists a character $\chi$
such that $\eta(\phi) \equiv \eta_{\chi}(\phi)$ (mod $\lambda$)
for all $\phi \in Z(\mathcal{H}_{K,\Z})$. $\qedsymbol$

\subsection{Proof of the Main Theorem}

We can now prove the more precise version of Theorem 2 alluded to
in the introduction.

\begin{thm}
Let $K=\prod_{v < \infty}K_v \subset G(\A^{\infty})$ be a compact
open subgroup. Let $\lambda|\ell$ be a finite place of $\bar{\Q}$
such that there exists at least two finite places $v$ where $\ell
\nmid |K_v|$ (this is automatic if $\ell > [F:\Q]n+1$). Let
$\pi=\otimes \pi_v$ be an automorphic representation of $G(\A)$
such that $\pi_{\infty}={\bf{1}}$ and $\pi^K \neq 0$. Assume $\pi$
is non-abelian modulo $\lambda$ relative to $K$. Let $w$ be a
finite place of $F$ such that $K_w$ is hyperspecial, and let
$J_w=K_w \cap K_w'$ be a parahoric subgroup, where $K_w' \neq K_w$
is maximal compact. Let $J=J_wK^w$ and $K'=K_w'K^w$. Suppose $\ell
\nmid [K':J]$ and
$$
\text{$\eta_{\pi^{K}}(e_{K,K'}) \equiv \eta_{{\bf{1}}}(e_{K,K'})$
(mod $\lambda$),}
$$
where
$$
e_{K,K'}=[K:J][K':J](e_K \star e_{K'} \star e_K) \in
Z(\mathcal{H}_{K,\Z}).
$$
Then there exists an automorphic representation
$\tilde{\pi}=\otimes \tilde{\pi}_v$ of $G(\A)$ such that
$\tilde{\pi}_{\infty}={\bf{1}}$ and $\tilde{\pi}^{K^w} \neq 0$
satisfying the following:

\begin{itemize}
\item $\tilde{\pi}_w^{J_w} \neq \tilde{\pi}_w^{K_w}+\tilde{\pi}_w^{K_w'}$, and \\
\item $\eta_{\tilde{\pi}^J}(\phi) \equiv \eta_{\pi^K}(e_K \star
\phi)$ (mod $\lambda$), for all $\phi \in Z(\mathcal{H}_{J,\Z})$.
\end{itemize}

\end{thm}

$Proof$. The reduction $\bar{\eta}_{\pi^K}$ modulo $\lambda \cap
\mathcal{O}$ factors through $\T_{J,\Z}^{\text{new}}$ by the main
lemma (Lemma 7). That is, there exists a character
$\eta':Z(\mathcal{H}_{J,\Z}) \rightarrow \F$ factoring through
$\T_{J,\Z}^{\text{new}}$ such that $\eta'(\phi)=\eta_{\pi^K}(e_K
\star \phi)$ (mod $\lambda$) for all $\phi \in
Z(\mathcal{H}_{J,\Z})$. As above, there is a surjective
homomorphism with nilpotent kernel
$$
\F_{\ell} \otimes_{\Z} \T_{J,\Z}^{\text{new}} \twoheadrightarrow
\T_{J,\F_{\ell}}^{\text{new}}.
$$
Thus $\eta'$ gives rise to a character $\T_{J,\F}^{\text{new}}
\rightarrow \F$, also denoted by $\eta'$. By a standard argument
(used above in section 6.2), there is an eigenform $f' \in
\mathcal{A}_{J,\F}^{\text{new}}$ with character $\eta'$. Now we
apply the Deligne-Serre lifting lemma, [DS, p. 522], to the finite
free module $\mathcal{A}_{J,\mathcal{O}_{\lambda}}^{\text{new}}$
where $\mathcal{O}_{\lambda}$ is the completion of $\mathcal{O}$
at $\lambda \cap \mathcal{O}$. It gives the existence of a
character $\tilde{\eta}: \T_{J,\Z}^{\text{new}} \rightarrow
\tilde{\mathcal{O}}_{\lambda}$ occurring in
$\mathcal{A}_{J,\tilde{\mathcal{O}}_{\lambda}}^{\text{new}}$ and
reducing to $\eta'$, where $\tilde{\mathcal{O}}_{\lambda}$ is the
ring of integers in a finite extension of $L_{\lambda}$. Since
$\T_{J,\Z}^{\text{new}}$ preserves the lattice
$\mathcal{A}_{J,\Z}^{\text{new}}$, the values $\tilde{\eta}(\phi)$
all lie in the ring of integers of some number field,
$\mathcal{O}_{\tilde{L}}$. We deduce that there exists a character
$\tilde{\eta}: Z(\mathcal{H}_{J,\Z}) \rightarrow
\mathcal{O}_{\tilde{L}}$, occurring in
$\mathcal{A}_{J}^{\text{new}}$, such that
$$
\text{$\tilde{\eta}(\phi) \equiv \eta_{\pi^K}(e_K \star \phi)$
(mod $\lambda$)}
$$
for all $\phi \in Z(\mathcal{H}_{J,\Z})$. From the decomposition
of $\mathcal{A}_J$ in terms of automorphic representations, it
follows that the newspace $\mathcal{A}_{J}^{\text{new}}$ has the
following description:
$$
\mathcal{A}_J^{\text{new}} \simeq {\bigoplus}_{\pi \in
\Pi_{\text{unit}}(G(\A)) : \pi_{\infty}={\bf{1}}}
m(\pi)(\pi^J/\pi^K+\pi^{K'}),
$$
as $Z(\mathcal{H}_J)$-modules. The center $Z(\mathcal{H}_J)$ acts
on the quotient $\pi^J/\pi^K+\pi^{K'}$ by the character
$\eta_{\pi^J}$. We conclude that there exists an automorphic
representation $\tilde{\pi}$ of $G(\A)$ with
$\tilde{\pi}_{\infty}={\bf{1}}$ and $\tilde{\pi}^J \neq
\tilde{\pi}^K+\tilde{\pi}^{K'}$, such that
$\eta_{\tilde{\pi}^J}=\tilde{\eta}$. In particular,
$$
\text{$\eta_{\tilde{\pi}^J}(\phi) \equiv \eta_{\pi^K}(e_K \star
\phi)$ (mod $\lambda$)}
$$
for all $\phi \in Z(\mathcal{H}_{J,\Z})$. This finishes the proof
$\qedsymbol$

\section{Applications}

\subsection{The Rank One Situation}

When the $F_w$-rank of $G_w^{\text{der}}$ is one, the condition
$\tilde{\pi}_w^{J_w} \neq
\tilde{\pi}_w^{K_w}+\tilde{\pi}_w^{K_w'}$ forces $\tilde{\pi}_w$
to be ramified:

\begin{cor}
With notation as above, let $w$ be a finite place of $F$ such that
$K_w$ is hyperspecial and the $F_w$-rank of $G_w^{\text{der}}$ is
one. Let $I_w=K_w \cap K_w'$ be an Iwahori subgroup, where $K_w'
\neq K_w$ is maximal compact. Let $I=I_wK^w$ and $K'=K_w'K^w$.
Suppose $\ell$ does not divide $[K:I][K':I]$ and
$$
\text{$\eta_{\pi^{K}}(e_{K,K'}) \equiv \eta_{{\bf{1}}}(e_{K,K'})$
(mod $\lambda$),}
$$
with $e_{K,K'}$ as in theorem 5. Then there exists an automorphic
representation $\tilde{\pi}=\otimes \tilde{\pi}_v$ of $G(\A)$ such
that $\tilde{\pi}_{\infty}={\bf{1}}$ and $\tilde{\pi}^{K^w} \neq
0$ satisfying the following:

\begin{itemize}
\item $\tilde{\pi}_w^{I_w} \neq 0$ and $\tilde{\pi}_w^{K_w}=0$, \\
\item $\eta_{\tilde{\pi}^I}(\phi) \equiv \eta_{\pi^K}(e_K \star
\phi)$ (mod $\lambda$), for all $\phi \in Z(\mathcal{H}_{I,\Z})$.
\end{itemize}

\end{cor}

$Proof$. Let $\tilde{\pi}$ be the automorphic representation we
get from the main theorem. We need to show $\tilde{\pi}_w$ is
ramified, so suppose on the contrary that $\tilde{\pi}_w^{K_w}
\neq 0$. Then $\tilde{\pi}_w^{K_w'} \neq 0$: The action of
$e_{K,K'}$ on $\tilde{\pi}^K$ factors through $\tilde{\pi}^{K'}$,
so if this is zero $\ell$ must divide $[K:I][K':I]$. Now, since
$\tilde{\pi}_w$ is a constituent of an unramified principal series
$\dim \tilde{\pi}_w^{I_w}$ is bounded by $|W|=2$. Consequently,
$\tilde{\pi}_w^{K_w} \cap \tilde{\pi}_w^{K_w'} \neq 0$. That is,
$\tilde{\pi}_w$ has nonzero $G_w^0$-invariants, and hence $\dim
\tilde{\pi}_w=1$. $\tilde{\pi}$ is automorphic, so by the strong
approximation theorem it must be one-dimensional. However, $\pi
\equiv \tilde{\pi}$ is assumed not to be abelian. $\qedsymbol$

$\e$

This corollary is a slight generalization of Bellaiche's theorem
1.4.6, [Bel, p. 215]: It gives results modulo arbitrary
$\lambda|\ell$, the level-raising condition is weaker, and we get
information about the action of the center of the Iwahori-Hecke
algebra on $\tilde{\pi}_w^{I_w}$. Bellaiche's proof is different.
He uses results of Lazarus and Vigneras from modular
representation theory, such as the computation of the composition
series of universal modules. With his stronger level-raising
condition, $\eta_{\pi^K}(\phi) \equiv \eta_{\bf{1}}(\phi)$ for
$all$ $\phi \in \mathcal{H}_{K_w}$, one can conclude that
$\tilde{\pi}_w$ is the actual Steinberg representation of $G_w$,
see [Bel, p. 221].

\subsection{$\U(3)$ - the Split Case}

In this subsection, we let $E/\Q$ denote an imaginary quadratic
extension of $\Q$, even though much of what we have to say is true
for CM extensions. We consider the quasi-split unitary $\Q$-group
in $3$ variables, $G^*=\U(2,1)$, split over $E$. We let $G=\U(3)$
be an arbitrary inner form of $G^*$ such that $G_{\infty}$ is
compact. Such exist since $E$ is imaginary. The rank is odd, so we
may even assume $G$ is quasi-split at all finite primes, but we do
not need that here. Now, we will focus on primes $q$ split in $E$.
First, we make some remarks on the parahoric subgroups of
$\GL_3(E_{\q}) \simeq \GL_3(\Q_q)$. There is the hyperspecial
maximal compact subgroup $K=\GL_3(\Z_q)$, and the Iwahori subgroup
$$
I=\{g \in K: \text{$g \equiv \begin{pmatrix} * & * & * \\ 0 & * &
* \\ 0 & 0 & * \end{pmatrix}$ (mod $q$) }\}.
$$
There is only one $\GL_3(\Q_q)$-conjugacy class of maximal proper
parahorics. We take
$$
\text{$J=\{g \in K: \text{$g \equiv \begin{pmatrix} * & * & * \\ *
&
* &
* \\ 0 & 0 & * \end{pmatrix}$ (mod $q$) }\}=K \cap \mu^{-1}K\mu$, $\y$ where
$\mu=\begin{pmatrix} q &  &  \\  & q &
 \\  &  & 1 \end{pmatrix}$,}
$$
as a representative. The following is a slightly stronger version
of Theorem 3.

\begin{thm}
Let $\pi=\otimes \pi_p$ be an automorphic representation of
$G(\A)$ with $\pi_{\infty}={\bf{1}}$. Let $\lambda|\ell$ be a
finite place of $\bar{\Q}$ such that $\pi$ is non-abelian modulo
$\lambda$. Choose a compact open subgroup $K=\prod K_p \subset
G(\A^{\infty})$ such that $\pi^K \neq 0$. If $\ell \leq 3$, or we
are in the situation where $E=\Q(\sqrt{-7})$ and $\ell=7$, assume
$\ell \nmid |K_p|$ for at least two primes $p$. Let $q \neq \ell$
be a prime, split in $E$, such that $K_q$ is hyperspecial. If
$\ell \nmid 1+q+q^2$, and the following is satisfied
$$
\text{$\h_{\pi,\q}\equiv
\begin{pmatrix}q & &
\\ & 1 & \\ & & q^{-1}
\end{pmatrix}$ (mod $\lambda$),}
$$
where $\q|q$, then there exists an automorphic representation
$\tilde{\pi}=\otimes \tilde{\pi}_p$ of $G(\A)$ with
$\tilde{\pi}_{\infty}={\bf{1}}$ and $\tilde{\pi}^{K^q}\neq 0$
satisfying the following conditions,

\begin{itemize}

\item $\tilde{\pi}_q$ is either an irreducible unramified
principal series or induced from a Steinberg representation. In
particular $\tilde{\pi}_q$ is generic, not $L^2$, and
$\tilde{\pi}_q^{J_q}\neq 0$, \\
\item $\eta_{\tilde{\pi}^J}(\phi)\equiv \eta_{\pi^K}(e_K \star
\phi)$ (mod $\lambda$), for all $\phi \in Z(\mathcal{H}_{J,\Z})$,
where $J=J_qK^q$.

\end{itemize}

\end{thm}

$Proof$. We first need to classify all the Iwahori-spherical
representations of $\GL_3(\Q_q)$. It is a theorem of Borel and
Casselman that these are precisely the constituents of unramified
principal series. Using the theory developed by Bernstein and
Zelevinsky, nicely summarized in [Kud], we obtain the following
table: $\nu=|\cdot|$ is the absolute value,

\center

\begin{tabular}{|c|c|c|c|c|c|c|c|}
\hline
   &  & constituent of & representation & unitary & tempered & $L^2$ & generic \\
 \hline \hline
  I &  & $\chi_1 \times \chi_2 \times \chi_3$ & $\chi_1 \times \chi_2 \times \chi_3$
  & below & $|\chi_i|=1$ &   & $\bullet$  \\
\hline
  II & a & $\chi_1\nu^{1/2}\times\chi_1\nu^{-1/2}\times \chi_2$
  & $\chi_1\text{St}_{\GL(2)}\times \chi_2$ & $|\chi_i|=1$  & $|\chi_i|=1$ &  & $\bullet$  \\
\cline{4-8}
   & b & $\chi_1\chi_2^{-1}\neq \nu^{\pm 3/2}$ & $\chi_1{\bf{1}}_{\GL(2)}\times \chi_2$
   & $|\chi_i|=1$  &   &   &   \\
 \hline
  III & a & $\chi\nu \times \chi \times \chi\nu^{-1}$ & $\chi\text{St}_{\GL(3)}$
  & $|\chi|=1$  & $|\chi|=1$ & $|\chi|=1$ & $\bullet$  \\
\cline{4-8}
  & b &  & $\chi V_P$ &  &  &  &  \\
 \cline{4-8}
 & c &  & $\chi V_Q$ &  &  &  &  \\
 \cline{4-8}
 & d &  & $\chi {\bf{1}}_{\GL(3)}$ & $|\chi|=1$ &  &  &  \\
 \hline
\end{tabular}

{\center{Table A: Iwahori-spherical representations of $\GL(3)$}}

\raggedright

$\et$

The irreducible representation $\chi_1 \times \chi_2 \times
\chi_2$ in group I is unitary if and only if either all the
$\chi_i$ are unitary, or, $\chi_1\chi_2^{-1}=\nu^{\alpha}$ with
$0<\alpha<1$ and $\chi_3$ unitary (after a permutation). In the
table, $P$ and $Q$ denote the parabolics of $G=\GL_3(\Q_q)$ of
type $(2,1)$ and $(1,2)$ respectively. Moreover, $V_P=C^{\infty}(P
\backslash G)/\C$ and $V_Q$ is defined similarly. They are not
unitary, and therefore irrelevant for the theory of automorphic
forms. Next, we list the dimensions of their parahoric fixed
spaces:

\center

\begin{tabular}{|c|c|c|c|c|c|c|}
\hline
   &   & representation & remarks & $K$ & $J$ & $I$ \\
 \hline \hline
  I &   & $\chi_1 \times \chi_2 \times \chi_3$
  &  & 1 &  3 & 6  \\
\hline
  II & a
  & $\chi_1\text{St}_{\GL(2)}\times \chi_2$ &   & 0 & 1 & 3  \\
\cline{3-7}
   & b & $\chi_1{\bf{1}}_{\GL(2)}\times \chi_2$
   &   &  1 &  2 & 3  \\
 \hline
  III & a & $\chi\text{St}_{\GL(3)}$
  &   & 0 & 0 & 1  \\
\cline{3-7}
  & b & $\chi V_P$ & not unitary & 0 & 1 & 2 \\
 \cline{3-7}
 & c & $\chi V_Q$ & not unitary & 0 & 1 & 2 \\
 \cline{3-7}
 & d & $\chi{\bf{1}}_{\GL(3)}$  & irrelevant & 1 & 1 & 1 \\
 \hline
\end{tabular}

{\center{Table B: Dimensions of the parahoric fixed spaces}}

\raggedright

$\et$

To compute these dimensions, we use the following observation: If
$P$ is parabolic and $J$ is parahoric, a choice of representatives
$g \in P \backslash G / J$ determines an isomorphism
$$
\Ind_P^G(\tau)^J \simeq {\bigoplus}_{g \in P \backslash G / J}
\tau^{P \cap gJg^{-1}},
$$
for every representation $\tau$ of a Levi factor $M_P$. In
particular, if $P=B$ is the Borel subgroup and $\tau$ is an
unramified character, the dimension of $\Ind_B^G(\tau)^J$ equals
the number of double cosets $|B \backslash G / J|$. With this
information, the proof proceeds as follows: Our main theorem gives
us an automorphic representation $\tilde{\pi}$ congruent to $\pi$
(modulo $\lambda$) such that $\tilde{\pi}_q^{J_q} \neq
\tilde{\pi}_q^{K_q}+\tilde{\pi}_q^{K_q'}$. Since $\tilde{\pi}_q$
must be unitary, we see from table B that it is of type I or IIa.
Then, from table A, we derive that $\tilde{\pi}_q$ is generic and
not $L^2$. Finally, note that there is a bijection $K/J \simeq
\GL_3(\F_q)/\bar{P}$, so $[K:J]=1+q+q^2$. $\qedsymbol$

$\e$

$Remark$. This corollary has no content unless $\pi_q$ is induced
from the determinant (type IIb), that is, unramified and
non-generic (and not $1$-dimensional), which is the case for the
endoscopic lifts from $\U(2) \times \U(1)$ considered in [Bel, p.
250]. In fact, the results we get for $\U(n)$ indicate that an
endoscopic abelian lift $\pi$ is congruent to a $\tilde{\pi}$
which is not endoscopic abelian. In his thesis [Bel, p. 218],
Bellaiche also has a result in the split case. Apparently, if you
only allow $\ell$ outside a finite set and $\pi$ occurs with
multiplicity $1$, then you can obtain a $\tilde{\pi}$ with
$\tilde{\pi}_q$ ramified. Hence, from our analysis,
$\tilde{\pi}_q$ is induced from Steinberg. It looks like the
preceding corollary is related to the $n=3$ case of conjecture 5.3
in [Ta2, p. 35], providing an analogue of Ihara's lemma, and to
the work of Mann [Man]. We also note that automorphic
representations of unitary groups with a generic component at a
split prime, come up naturally in the proof of the local Langlands
correspondence for $\GL(n)$ [HT].

\subsection{$\GSp(4)$}

In this subsection we view $\GSp(4)$ as an algebraic $\Q$-subgroup
of $\GL(4)$ by realizing it with respect to the standard
skew-diagonal symplectic form. With this choice, the set of upper
triangular matrices form a Borel subgroup $B=TU$. There are two
maximal parabolic subgroups containing $B$, namely the Siegel
parabolic
$$
P=M_P \ltimes N_P=\{\begin{pmatrix} g & \\ & \nu {^{\tau}g^{-1}}
\end{pmatrix}
\begin{pmatrix}1 & & r & s \\ & 1 & t & r \\ & & 1 & \\ & & & 1
\end{pmatrix}\},
$$
where ${^{\tau}g}$ denotes the skew-transpose, and the Klingen
parabolic
$$
Q=M_Q \ltimes N_Q=\{\begin{pmatrix} \nu & &   \\ & g &   \\ & &
\nu^{-1}\det g\end{pmatrix}
\begin{pmatrix} 1& c& &  \\ &1 &  & \\ &  & 1&-c \\  & & & 1
\end{pmatrix}
\begin{pmatrix} 1& & r& s \\ &1 &  & r\\ &  & 1& \\  & & & 1
\end{pmatrix}  \}.
$$
We consider an inner form $G$ of $\GSp(4)$ such that
$G^{\text{der}}(\R)$ is compact. Concretely we have $G=\GSpin(f)$,
where $f$ is some definite quadratic form in $5$ variables over
$\Q$. Now, let us first describe the parahoric subgroups of
$\GSp_4(\Q_q)$. There is the hyperspecial maximal compact subgroup
$K_q=\GSp_4(\Z_q)$, and the Iwahori subgroup $I_q$ consisting of
elements in $K_q$ with upper triangular reduction mod $q$.
Similarly, $P$ and $Q$ define (non-conjugate) parahoric subgroups
$J_q'$ and $J_q$ called the Siegel parahoric and the Klingen
parahoric respectively. One can easily check that we have the
identity,
$$
\text{$J_q'=K_q \cap hK_qh^{-1}$, where $h=\begin{pmatrix} & & 1&  \\
& & &1
\\ q&  & & \\  & q& &
\end{pmatrix}$.}
$$
However, $J_q=K_q \cap K_q'$, where $K_q'$ is the non-special
maximal compact subgroup containing $I_q$. It is called the
paramodular group. Since $P$ and $Q$ are not associated
parabolics, the classification of the Iwahori-spherical
representations of $\GSp_4(\Q_q)$ is much more complicated than
for $\GL_3(\Q_q)$. Fortunately, this has been done by Ralf
Schmidt. The tables we need are Table 1 and Table 3 in the
forthcoming paper [Sch]. With the permission of Ralf Schmidt, we
have reproduced the information we need in Appendix 2 as Table C
and Table D. We use the notation from this appendix below. If
$\pi$ has a Galois representation $\rho_{\pi,\lambda}$ (for
example, if it transfers to a cuspidal representation $\Pi$ of
$\GSp(4)$ with $\Pi_{\infty}$ in the discrete series, see [Lau],
[Wei]), then $\rho_{\pi,\lambda}(\Fr_p)$ and $\h_{\pi_p\otimes
|\nu|^{-3/2}}$ have the same eigenvalues. In this case, $\pi$ is
abelian modulo $\lambda$ if some twist of
$\bar{\rho}_{\pi,\lambda}$ has the form
${\bf{1}}\oplus\bar{\omega}_{\ell}\oplus\bar{\omega}_{\ell}^2
\oplus \bar{\omega}_{\ell}^3$. We obtain the following
strengthening of Theorem 4.

\begin{thm}
Let $\pi=\otimes \pi_p$ be an automorphic representation of
$G(\A)$ with $\pi_{\infty}={\bf{1}}$. Let $\lambda|\ell$ be a
finite place of $\bar{\Q}$ such that $\pi$ is non-abelian modulo
$\lambda$. Choose a compact open subgroup $K=\prod K_p$ such that
$\pi^K \neq 0$. If $\ell \leq 5$ assume $\ell \nmid |K_p|$ for at
least two primes $p$. Let $q \neq \ell$ be a prime such that $K_q$
is hyperspecial. Suppose
$$
\text{$\h_{\pi_q\otimes |\nu|^{-3/2}}\equiv
\begin{pmatrix}1 & & &
\\ & q &  & \\ & & q^2 & \\ & & & q^3
\end{pmatrix}$ (mod $\lambda$).}
$$
Then there exists an automorphic representation
$\tilde{\pi}=\otimes \tilde{\pi}_p$ of $G(\A)$ with
$\tilde{\pi}_{\infty}={\bf{1}}$ and $\tilde{\pi}^{K^q}\neq 0$
satisfying the following conditions,

\begin{itemize}

\item $\tilde{\pi}_q$ is generic and Klingen-spherical,\\
\item $\eta_{\tilde{\pi}^J}(\phi)\equiv \eta_{\pi^K}(e_K \star
\phi)$ (mod $\lambda$), for all $\phi \in Z(\mathcal{H}_{J,\Z})$,
where $J=J_qK^q$.

\end{itemize}

Moreover, if in addition $q^4 \neq 1$ (mod $\ell$),
$\tilde{\pi}_q$ must be of type I, IIa or IIIa.
\end{thm}

$Proof$. We apply the main theorem (Theorem 5) to the Klingen
parahoric $J_q$. An easy computation shows that $[K_q':J_q]=q$. We
get an automorphic representation $\tilde{\pi}$, congruent to
$\pi$ modulo $\lambda$, such that the component at $q$ satisfies
the identity:
$$
\tilde{\pi}_q^{J_q}\neq \tilde{\pi}_q^{K_q}+\tilde{\pi}_q^{K_q'}.
$$
In particular, $\tilde{\pi}_q^{J_q}\neq 0$. We must have that
$\tilde{\pi}_q^{K_q} \cap \tilde{\pi}_q^{K_q'}=0$, for otherwise
$\dim \tilde{\pi}_q=1$ and therefore $\tilde{\pi}$ is
one-dimensional by the strong approximation theorem. However,
$\pi$ is assumed to be non-abelian modulo $\lambda$. Thus,
equivalently we have
$$
\dim \tilde{\pi}_q^{J_q} > \dim \tilde{\pi}_q^{K_q}+\dim
\tilde{\pi}_q^{K_q'}.
$$
From Schmidt's tables, [Sch, p. 16] (that is, Table D in Appendix
2), we deduce that this inequality is satisfied precisely when
$\tilde{\pi}_q$ is of type I, IIa, IIIa, IVb, IVc, Va or VIa.
However, those representations of type IVb and IVc are not unitary
and can therefore be ruled out immediately. We are then left with
the possible types I, IIa, IIIa, Va and VIa. Then, from the tables
[Sch, p. 9] (Table C in Appendix 2), we read off that
$\tilde{\pi}_q$ is generic. Indeed all the representations of type
Xa are generic, for X arbitrary.

$\e$

Now, let us show that the types Va and VIa can also be ruled out
if we assume $q^4 \neq 1$ (mod $\ell$). Suppose first that
$\tilde{\pi}_q$ is of type Va, that is, the unique
subrepresentation of some $|\cdot|\xi_0 \times \xi_0 \rtimes
|\cdot|^{-1/2}\sigma$ where $\xi_0$ has order two, see [Sch, p. 7]
for an explanation of the notation. By the main theorem, the
center of the Klingen-Hecke algebra $Z(\mathcal{H}_{J_q,\Z})$ acts
on $\tilde{\pi}_q^{J_q}$ by a character
$\eta_{\tilde{\pi}_q^{J_q}}$ satisfying the congruence
$$
\text{$\eta_{\tilde{\pi}_q^{J_q}}(\phi) \equiv
\eta_{\pi_q^{K_q}}(e_{K_q}\star \phi) \equiv
\eta_{{\bf{1}}}(e_{K_q}\star \phi)$ (mod $\lambda$)},
$$
for all $\phi \in Z(\mathcal{H}_{J_q,\Z})$. We get immediately
that the analogous statement is also true for the center of the
Iwahori-Hecke algebra $Z(\mathcal{H}_{I_q,\Z})$. This, however,
acts by a character on the Iwahori-fixed vectors in the principal
series $|\cdot|\xi_0 \times \xi_0 \rtimes |\cdot|^{-1/2}\sigma$
(for it has an unramified Langlands quotient, so is generated by
any nonzero $K_q$-fixed vector). Hence, $Z(\mathcal{H}_{I_q,\Z})$
acts on every constituent of this principal series by the same
character $\eta_{\tilde{\pi}_q^{I_q}}$. In particular, the action
of the spherical Hecke algebra $\mathcal{H}_{K_q,\Z} \simeq
Z(\mathcal{H}_{I_q,\Z})$ on the $K_q$-fixed vectors of the
unramified quotient (type Vd) is given by a character congruent to
$\eta_{{\bf{1}}}$. In terms of their Satake parameters we
therefore must have (modulo the action of the Weyl group):
$$
\text{$\begin{pmatrix} q^{-1/2}\sigma(q) & & & \\
& q^{-1/2}\xi_0\sigma(q) & & \\
& & q^{1/2}\xi_0\sigma(q) & \\
& & & q^{1/2}\sigma(q)\end{pmatrix}\equiv
\begin{pmatrix} q^{-3/2} & & & \\
& q^{-1/2} & & \\
& & q^{1/2} & \\
& & & q^{3/2}\end{pmatrix}$ (mod $\lambda$).}
$$
Since $\xi_0(q)=-1$ we conclude that $q\equiv -1$ or $q^2 \equiv
-1$ modulo $\ell$. Secondly, assume $\tilde{\pi}_q$ is of type
VIa, that is, the unique irreducible subrepresentation of some
$|\cdot| \times {\bf{1}} \rtimes |\cdot|^{-1/2}\sigma$. Then, by
the argument above, we conclude that the unramified quotient of
this principal series must be congruent to ${\bf{1}}$. That is, in
terms of their Satake parameters:
$$
\text{$\begin{pmatrix} q^{-1/2}\sigma(q) & & & \\
& q^{-1/2}\sigma(q) & & \\
& & q^{1/2}\sigma(q) & \\
& & & q^{1/2}\sigma(q)\end{pmatrix}\equiv
\begin{pmatrix} q^{-3/2} & & & \\
& q^{-1/2} & & \\
& & q^{1/2} & \\
& & & q^{3/2}\end{pmatrix}$ (mod $\lambda$).}
$$
It follows that $q^2 \equiv 1$. The types I, IIa and IIIa cannot
be excluded, even if $\pi$ has trivial central character.
$\qedsymbol$

$\e$

$Remark$. There exists $q$ with $q^4 \neq 1$ (mod $\ell$)
precisely when $\ell \geq 7$. In this case $\tilde{\pi}_q$ is an
unramified principal series (type I) or induced from a twisted
Steinberg representation $\chi \text{St}_{\GL(2)} \rtimes \chi'$
or $\chi \rtimes \chi' \text{St}_{\GL(2)}$ (type IIa and IIIa
respectively). If one can show that $\tilde{\pi}_q$ is
para-ramified, meaning that $\tilde{\pi}_q$ has no nonzero
$K_q'$-fixed vectors, one can conclude that it is of type IIIa and
therefore induced from a twisted Steinberg representation on the
Klingen-Levi. It seems possible to prove this if $m(\pi)=1$, using
the methods of [Bel] and [Clo]. We hope to return to this point in
another paper. The result above only gives non-trivial congruences
if $\pi_q$ is non-generic. If $\pi$ is of Saito-Kurokawa type
(that is, a theta-lift from the $\widetilde{\SL}(2)$), it is
locally non-generic, and we get a $\tilde{\pi}$ congruent to $\pi$
which is not of Saito-Kurokawa type. If we know $\tilde{\pi}_q$ is
of type IIIa, we can apply this strategy to the Bloch-Kato
conjecture for the motives attached to classical modular forms of
weight (at least) $4$, using the methods of [Bel]. We should note
that if we choose to work with the Siegel-parahoric $J_q'$, we can
only conclude that $\tilde{\pi}_q$ is generic $or$ a
Saito-Kurokawa lift.

\section*{Appendix 1. Congruent Representations}

The compact open subgroups $K \subset G(\A^{\infty})$ form a
directed set by opposite inclusion, that is $K \preccurlyeq J
\Leftrightarrow K \supset J$. Let $R$ be a commutative ring. As
$K$ varies over the compact open subgroups, the centers
$Z(\mathcal{H}_{K,R})$ form an inverse system of $R$-algebras with
respect to the canonical maps $Z(\mathcal{H}_{K,R}) \leftarrow
Z(\mathcal{H}_{J,R})$ when $K \supset J$. Let
$$
\mathcal{Z}_{G(\A^{\infty}),R}=\underleftarrow \lim
Z(\mathcal{H}_{K,R}).
$$
In this limit, it is enough to let $K$ run through a neighborhood
basis at the identity. Thus $\mathcal{Z}_{G(\A^{\infty}),R}$ is a
commutative $R$-algebra, and it comes with projections ($K \supset
J$)

\begin{align*}
\xymatrix{  & \mathcal{Z}_{G(\A^{\infty}),R} \ar[dl]_{\text{pr}_K}
\ar[dr]^{\text{pr}_J}
& \\
 Z(\mathcal{H}_{K,R}) & &  Z(\mathcal{H}_{J,R})\ar[ll]^{e_K \star \phi \leftarrow \phi}}
\end{align*}

All we have said makes sense for any locally profinite group, so
in particular we have local analogues $\mathcal{Z}_{G_v,R}$ for
each finite place $v$. If $\mu=\otimes \mu_v$, it follows that
$$
\mathcal{Z}_{G(\A^{\infty}),R} \simeq \bigotimes_{v <
\infty}\mathcal{Z}_{G_v,R},
$$
a restricted tensor product. Indeed the decomposable groups
$K=\prod K_v$ form a cofinal system. It remains to determine the
algebras $\mathcal{Z}_{G_v,R}$. By [Cas, p. 14], there exists a
neighborhood basis at $1$ consisting of compact open subgroups
$K_v \subset G_v$ with Iwahori factorization with respect to a
fixed minimal parabolic. If $G_v$ is unramified, for such a $K_v$
the canonical map $Z(\mathcal{H}_{K_v,R})\rightarrow
\mathcal{H}_{v,R}^{\text{sph}}$ to the spherical Hecke algebra at
$v$ is an isomorphism [Bu1], [Bu2]. This is a well-known result
due to Bernstein when $K_v$ is an actual Iwahori subgroup.
Therefore,
$$
\text{$G_v$ unramified $\Longrightarrow$ $\mathcal{Z}_{G_v,R}
\simeq \mathcal{H}_{v,R}^{\text{sph}}$.}
$$
The reason for introducing these objects is the following: Let
$\pi=\otimes \pi_v$ be an irreducible admissible representation of
$G(\A)$. Then there exists a unique character
$$
\eta_{\pi}:\mathcal{Z}_{G(\A^{\infty}),\Z} \rightarrow \C,
$$
such that $\eta_{\pi}=\eta_{\pi^K}\circ \text{pr}_K$ for every $K$
such that $\pi^K \neq 0$. Uniqueness is clear, and the existence
reduces to showing that $\eta_{\pi^J}(\phi)=\eta_{\pi^K}(e_K \star
\phi)$ for $K \supset J$ when $\pi^K \neq 0$. Similarly, we have
characters $\eta_{\pi_v}$ locally, and $\eta_{\pi}=\otimes
\eta_{\pi_v}$ under the isomorphism above. If $\pi$ is automorphic
and $\pi_{\infty}={\bf{1}}$, the character $\eta_{\pi}$ maps into
the ring of integers of some number field. Our work suggests the
following definition:

\begin{df}
Let $\pi$ and $\tilde{\pi}$ be automorphic representations of
$G(\A)$, both trivial at infinity, and let $\lambda$ be a finite
place of $\bar{\Q}$. Then we say that $\pi$ and $\tilde{\pi}$ are
congruent modulo $\lambda$, and we write $\tilde{\pi}\equiv \pi$
(mod $\lambda$), if for all $\phi \in
\mathcal{Z}_{G(\A^{\infty}),\Z}$ we have
$$
\text{$\eta_{\tilde{\pi}}(\phi)\equiv \eta_{\pi}(\phi)$ (mod
$\lambda$).}
$$
\end{df}

Analogously, it makes sense to say the local components
$\tilde{\pi}_v$ and $\pi_v$ are congruent. Then $\tilde{\pi}\equiv
\pi$ (mod $\lambda$) if and only $\tilde{\pi}_v\equiv \pi_v$ (mod
$\lambda$) for all $v < \infty$. This is the kind of local-global
compatibility aimed for in Parson's thesis [Par]. Parson has
another definition of being congruent. We do not know how the two
definitions are related. Note also that if $\tilde{\pi}_v$ and
$\pi_v$ are both unramified, then $\tilde{\pi}_v\equiv \pi_v$ (mod
$\lambda$) means that the Satake parameters are congruent as it
should. With these definitions, our results translate into those
stated in the introduction.

\section*{Appendix 2. Iwahori-Spherical Representations of $\GSp(4)$}

In this appendix we reproduce parts of Table 1 and Table 3 in
[Sch]. We are grateful to Ralf Schmidt for his permission to do
so. We stress that the tables in [Sch] contain more information
than what is listed here (such as Atkin-Lehner eigenvalues and
signs of $\epsilon$-factors). Below, we employ the notation of
[ST]. Thus $\nu$ denotes the normalized absolute value of a
non-archimedean local field. If $\chi_1$, $\chi_2$ and $\sigma$
are unramified characters, we recall that $\chi_1 \times \chi_2
\rtimes \sigma$ denotes the principal series of $\GSp(4)$ obtained
from
$$
T \ni \diag(x,y,zy^{-1},zx^{-1}) \mapsto
\chi_1(x)\chi_2(y)\sigma(z) \in \C^*
$$
by normalized induction. Similarly, if $\pi$ is a representation
of $\GL(2)$, we denote by $\pi \rtimes \sigma$ and $\sigma \rtimes
\pi$ the representations of $\GSp(4)$ induced from
$\diag(X,z{^{\tau}X}^{-1}) \mapsto \pi(X)\sigma(z)$ and
$\diag(z,X,z^{-1}\det X) \mapsto \sigma(z)\pi(X)$ respectively. By
$L((-))$ we mean the unique irreducible quotient (the Langlands
quotient) when it exists. The representations
$\tau(S,\nu^{-1/2}\sigma)$ and $\tau(T,\nu^{-1/2}\sigma)$ are the
constituents of ${\bf{1}}\rtimes \sigma \text{St}_{\GL(2)}$. They
are occasionally called limit of discrete series. $\xi_0$ is the
non-trivial unramified quadratic character.

\center

\begin{tabular}{|c|c|c|c|c|c|c|}
\hline
   &  & constituent of & representation &  tempered & $L^2$ & generic \\
 \hline \hline
  I &  & $\chi_1 \times \chi_2 \rtimes \sigma$
  & $\chi_1 \times \chi_2 \rtimes \sigma$  & $|\chi_i|=|\sigma|=1$  &   & $\bullet$\\
\hline
  II & a & $\nu^{1/2}\chi \times \nu^{-1/2}\chi \rtimes \sigma$,
& $\chi\text{St}_{\GL(2)}\rtimes \sigma$  &  $|\chi|=|\sigma|=1$ &  & $\bullet$\\
\cline{4-7}
   & b & $\chi^2 \notin \{\nu^{\pm 1}, \nu^{\pm 3}\}$ & $\chi{\bf{1}}_{\GL(2)}\rtimes \sigma$
     &   &   &   \\
 \hline
  III & a & $\chi\times \nu \rtimes \nu^{-1/2}\sigma$, & $\chi \rtimes \sigma \text{St}_{\GL(2)}$
   &  $|\chi|=|\sigma|=1$ &   &   $\bullet$\\
\cline{4-7}
  & b & $\chi \notin \{{\bf{1}},\nu^{\pm 2}\}$ & $\chi \rtimes \sigma {\bf{1}}_{\GL(2)}$  &  &  &  \\
 \hline
 IV & a & $\nu^2 \times \nu \rtimes \nu^{-3/2}\sigma$ & $\sigma\text{St}_{\GSp(4)}$  & $\bullet$ & $\bullet$ & $\bullet$ \\
 \cline{4-7}
 & b &  & $L((\nu^2,\nu^{-1}\sigma\text{St}_{\GL(2)}))$   &  &  &  \\
\cline{4-7}
 & c &  & $L((\nu^{3/2}\text{St}_{\GL(2)},\nu^{-3/2}\sigma))$  &  &  &  \\
\cline{4-7}
 & d &  & $\sigma{\bf{1}}_{\GSp(4)}$  &  &  &  \\
\hline
V & a & $\nu\xi_0 \times \xi_0 \rtimes \nu^{-1/2}\sigma$, & $\delta([\xi_0,\nu\xi_0],\nu^{-1/2}\sigma)$  & $\bullet$ & $\bullet$ &  $\bullet$\\
 \cline{4-7}
 & b & $\xi_0^2={\bf{1}}$, $\xi_0 \neq {\bf{1}}$ & $L((\nu^{1/2}\xi_0\text{St}_{\GL(2)},\nu^{-1/2}\sigma))$  &  &  &  \\
\cline{4-7}
 & c &  & $L((\nu^{1/2}\xi_0\text{St}_{\GL(2)},\xi_0\nu^{-1/2}\sigma))$  &  &  &  \\
\cline{4-7}
 & d &  & $L((\nu\xi_0,\xi_0 \rtimes \nu^{-1/2}\sigma))$  &  &  &  \\
\hline
VI & a & $\nu \times {\bf{1}} \rtimes \nu^{-1/2}\sigma$ & $\tau(S,\nu^{-1/2}\sigma)$  & $\bullet$ &  &  $\bullet$\\
 \cline{4-7}
 & b &  & $\tau(T,\nu^{-1/2}\sigma)$  & $\bullet$ &  &  \\
\cline{4-7}
 & c &  & $L((\nu^{1/2}\text{St}_{\GL(2)},\nu^{-1/2}\sigma))$  &  &  &  \\
\cline{4-7}
 & d &  & $L((\nu, {\bf{1}} \rtimes \nu^{-1/2}\sigma ))$  &  &  &  \\
\hline

\end{tabular}

{\center{Table C: Iwahori-spherical representations of $\GSp(4)$}}

\raggedright

$\et$

In the following table, our notation is different from [Sch].
Recall that in our setup $K$ is hyperspecial, $K'$ is paramodular,
$J$ is the Klingen parahoric, $J'$ the Siegel parahoric and $I$ is
the Iwahori subgroup of $\GSp(4)$.

\center

\begin{tabular}{|c|c|c|c|c|c|c|c|c|}
\hline
   &   & representation & remarks & $K$ & $K'$ & $J$ & $J'$ & $I$ \\
 \hline \hline
  I &   &
   $\chi_1 \times \chi_2 \rtimes \sigma$&  & 1 &  2 & 4 & 4 & 8\\
\hline
  II & a
  & $\chi\text{St}_{\GL(2)}\rtimes \sigma$ &   & 0 & 1 & 2 & 1& 4\\
\cline{3-9}
   & b & $\chi{\bf{1}}_{\GL(2)}\rtimes \sigma$
   &   & 1  &  1 & 2 & 3& 4\\
 \hline
  III & a & $\chi \rtimes \sigma \text{St}_{\GL(2)}$
  &   &  0 & 0 & 1 & 2& 4\\
\cline{3-9}
  & b & $\chi \rtimes \sigma {\bf{1}}_{\GL(2)}$ &  &  1& 2 & 3 & 2& 4\\
 \hline
IV & a & $\sigma\text{St}_{\GSp(4)}$ &  & 0 &  0& 0 &0 & 1\\
 \cline{3-9}
 & b &  $L((\nu^2,\nu^{-1}\sigma\text{St}_{\GL(2)}))$ & not unitary & 0 &  0& 1 & 2& 3\\
 \cline{3-9}
 & c & $L((\nu^{3/2}\text{St}_{\GL(2)},\nu^{-3/2}\sigma))$  & not unitary & 0 & 1 &  2&1 &3 \\
 \cline{3-9}
 & d & $\sigma{\bf{1}}_{\GSp(4)}$  & irrelevant & 1 & 1 & 1 & 1& 1\\
\hline
V & a & $\delta([\xi_0,\nu\xi_0],\nu^{-1/2}\sigma)$ &  & 0 & 0 & 1 &0 & 2\\
 \cline{3-9}
 & b & $L((\nu^{1/2}\xi_0\text{St}_{\GL(2)},\nu^{-1/2}\sigma))$  &  & 0 & 1 & 1 &1 & 2 \\
 \cline{3-9}
 & c &  $L((\nu^{1/2}\xi_0\text{St}_{\GL(2)},\xi_0\nu^{-1/2}\sigma))$ &  & 0 & 1 & 1 &1 & 2\\
 \cline{3-9}
 & d & $L((\nu\xi_0,\xi_0 \rtimes \nu^{-1/2}\sigma))$  &  & 1 & 0 & 1 & 2& 2\\
\hline
VI & a & $\tau(S,\nu^{-1/2}\sigma)$ &  & 0 & 0 &1  & 1& 3\\
 \cline{3-9}
 & b & $\tau(T,\nu^{-1/2}\sigma)$  &  & 0 & 0 & 0 & 1& 1\\
 \cline{3-9}
 & c & $L((\nu^{1/2}\text{St}_{\GL(2)},\nu^{-1/2}\sigma))$  &  & 0 & 1 & 1 & 0& 1\\
 \cline{3-9}
 & d & $L((\nu, {\bf{1}} \rtimes \nu^{-1/2}\sigma ))$  &  & 1 & 1 & 2 &2 & 3\\
\hline

\end{tabular}

{\center{Table D: Dimensions of the parahoric fixed spaces}}

\raggedright

$\et$

\newpage


$\et$

\begin{itemize}
\item[] Claus Mazanti Sorensen,\\
253-37 Caltech,\\
Pasadena, CA 91125, USA.\\
{\texttt{claus@caltech.edu}}
\end{itemize}

\end{document}